\documentclass[opre,sglanonrev]{Arxiv}
\usepackage{lmodern}
\usepackage{eqndefns-left} 
\RequirePackage{tgtermes}
\RequirePackage{newtxtext}
\RequirePackage{newtxmath}
\RequirePackage{bm}
\RequirePackage{endnotes}

\OneAndAHalfSpacedXII 

\usepackage{booktabs} 
\usepackage{subcaption}
\usepackage{amsmath,amssymb,graphicx,url,mathtools}
\usepackage{xcolor}
\usepackage{soul}
\usepackage{subfiles}
\usepackage{pgfplots}
\pgfplotsset{compat=1.18}
\usepackage{tikz}
\usepackage{calc,tikz,nicefrac,xifthen,graphicx,subcaption,url}
\usetikzlibrary{patterns}
\usetikzlibrary{math}
\usetikzlibrary{shapes.misc}
\usetikzlibrary{decorations.pathreplacing}
\usetikzlibrary {decorations.pathmorphing}
\usepackage{algorithm,algorithmic}

\newcommand{\ProximalPointMethodFull}{proximal point method}
\newcommand{\ProximalPointMethodFullCap}{Proximal Point Method}
\newcommand{\ProximalPointMethodAbb}{PPM}

\newcommand{\ProximalPathFull}{proximal path}
\newcommand{\ProximalPathFullCap}{Proximal Path}
\newcommand{\ProximalPathAbb}{PP}

\newcommand{\CentralPathFull}{central path}
\newcommand{\CentralPathFullCap}{Central Path}
\newcommand{\CentralPathAbb}{CP}

\newcommand{\eg}{\emph{e.g.}}
\newcommand{\ie}{\emph{i.e.}}

\usepackage{natbib}
 \bibpunct[, ]{(}{)}{,}{a}{}{,}%
 \def\bibfont{\small}%
\EquationsNumberedThrough    

\TheoremsNumberedThrough     
\ECRepeatTheorems  %

\MANUSCRIPTNO{}

\begin{document}



\RUNTITLE{}

\TITLE{Robust Paths: Geometry and Computation}

\ARTICLEAUTHORS{%
\AUTHOR{Hao Hao}
\AFF{Carnegie Mellon University, \EMAIL{haohao@andrew.cmu.edu}}
\AUTHOR{Peter Zhang}
\AFF{Carnegie Mellon University, 
\EMAIL{pyzhang@cmu.edu}}

} 

\ABSTRACT{
Applying robust optimization often requires selecting an appropriate uncertainty set---both in shape and size---a choice that directly affects the trade-off between average-case and worst-case performances. In practice, this calibration is usually done via trial-and-error: solving the robust optimization problem many times with different uncertainty set shapes and sizes, and examining their performance trade-off. This process is computationally expensive and ad hoc. In this work, we take a principled approach to study this issue for robust optimization problems with linear objective functions, convex feasible regions, and convex uncertainty sets. We introduce and study what we define as the \emph{robust path}: a set of robust solutions obtained by varying the uncertainty set's parameters. Our central geometric insight is that a robust path can be characterized as a Bregman projection of a curve (whose geometry is defined by the uncertainty set) onto the feasible region. This leads to a surprising discovery that the robust path can be approximated via the trajectories of standard optimization algorithms, such as the {\ProximalPointMethodFull}, of the deterministic counterpart problem. We give a sharp approximation error bound and show it depends on the geometry of the feasible region and the uncertainty set. We also illustrate two special cases where the approximation error is zero: the feasible region is polyhedrally monotone (\eg, a simplex feasible region under an ellipsoidal uncertainty set), or the feasible region and the uncertainty set follow a dual relationship. We demonstrate the practical impact of this approach in two settings: portfolio optimization and adversarial deep learning. The former numerically validates the zero approximation error under favorable conditions (feasible region is polyhedrally monotone); and when the technical conditions are violated, still retains a very small error. The latter case severely breaks the linear objective condition in our theory. But our solution technology still shows strong performance: $85\%$ reduction in computational time and near-Pareto efficiency in terms of average-case and worst-case performances.
}

\maketitle

\section{Introduction}\label{sec:Intro}
We study robust optimization problems with uncertain objective functions:
\begin{equation}
    \label{eq: RC}
    (\mathrm{RC}) \quad \min_{x\in \mathcal{X}} \max_{a \in \mathcal{U}} \langle a, x \rangle,
\end{equation}
where decision $x$ belongs to a closed, convex and nonempty feasible region $\mathcal{X} \subseteq \mathbb{R}^{n}$, $a$ is a vector of uncertain parameters and is only known to reside in an uncertainty set: $$\mathcal{U}= \{a_0 + \xi: \xi \in \Xi \subset \mathbb{R}^{n}\}.$$
Here $a_0$ is the nominal vector, $\xi$ is the uncertain perturbation assumed to be in a compact and convex set $\Xi$. 
When deploying $(\mathrm{RC})$ in the real world, the design of the uncertainty set $\Xi$ is crucial. To this end, we consider general uncertainty sets representable via gauge function \citep{ freund1987dual, Gauge_Opt, wei2025redefiningcoherentriskmeasures} constraints:
$$\Xi(r,\mathcal{V})=  \left\{\xi \in \mathbb{R}^{n}:    \|\xi\|_{\mathcal{V}} \leq r  \right\},$$
where the gauge function is defined as $\|v\|_{\mathcal{V}}=\inf \{t \geq 0: v\in t\mathcal{V} \}$, and the size and shape of $\Xi(r,\mathcal{V})$ can be flexibly adjusted via radius $r$ and gauge set $\mathcal{V}$ respectively. We assume $r\geq0$ and $\mathcal{V}$ is a compact and convex set with $0 \in \mathrm{int}(\mathcal{V})$.  For instance, selecting $\mathcal{V}$ to be an ellipsoid or unit $l_{p}$ norm balls recovers ellipsoidal or $l_{p}$ norm type uncertainty sets.

The calibration of the uncertainty set size and shape $(r,\mathcal{V})$ directly affects the performance of the resulting robust solution, yet it is a difficult task in practice. The design of $(r,\mathcal{V})$ controls the robustness and efficiency trade-off of the deployed robust solution. At one end of the spectrum, if one assumes $r=0$, \ie, $a$ takes its nominal value $a_{0}$ deterministically, $(\mathrm{RC})$ is reduced to the following deterministic optimization problem:
\begin{equation}
    \label{eq: P}
    (\mathrm{P}) \quad \min_{x\in \mathcal{X}}  \langle a_{0}, x \rangle.
\end{equation}
A resulting optimal solution $x_{\mathrm{E}}$ (E for efficiency) performs well if indeed $a=a_{0}$; however, there is no performance guarantee when the true realization of $a$ deviates from $a_{0}$. At the other end of the spectrum, a highly risk-averse decision maker may select a large $r=\overline{r}$, and solve $(\mathrm{RC})$ for a robust solution, which is robustly optimal under $a\in\mathcal{U}(\overline{r}, \mathcal{V})$, yet it may perform poorly under $a=a_{0}$. The uncertainty set shape $\mathcal{V}$ is another important design lever. $\mathcal{V}$ is typically designed to leverage the available information on the uncertainty while ensuring the computational tractability of its resulting (RC). 

Existing approaches calibrate $(r, \mathcal{V})$ based on probabilistic guarantees \citep{Bertsimas2004, bertsimas2021probabilistic, mohajerin2018data, Blanchet_1}. However, such approaches assume prior knowledge of the uncertainty distribution or observations on the uncertainty, which can be unavailable in practice. Even with distributional information or data on the uncertainty, it has been observed that the resulting robust solutions from this approach can be too conservative \citep{sim2021new}. Often, the practical approach is costly and ad hoc: solving the robust optimization problem multiple times, each under a different choice of $(r,\mathcal{V})$, before comparing the performance of the different robust solutions (which includes cross-validation among other statistical methods for hyperparameter tuning) \citep{Ben_Phi_divergence, sim2021new, mohajerin2018data}. To this end, ideally, decision makers need \emph{the entire set of robust solutions} under multiple shapes $\mathcal{V}$ and radii $r$.
\begin{definition}[Robust Path]\label{def:robust_path} \textit{The \emph{robust path} of $(\mathrm{RC})$ under $\mathcal{V}$ is defined as}
$$\mathcal{P}(\mathcal{V}) = \left\{ x_{\mathrm{R}}(r, \mathcal{V}) \in \argmin_{x\in \mathcal{X}} \max_{\xi \in \Xi(r,\mathcal{V})} \langle a_{0} + \xi, x \rangle : r\in [0,\infty) \right\}.$$ 
\end{definition}
The challenge lies in obtaining the robust paths, potentially without repeatedly solving the robust counterpart.  For any hope in tackling this problem, one first needs to study the structure of the robust paths under different $\mathcal{V}$, before exploiting the structural information to find an algorithm capable of tracing the robust paths. To this end, we investigate the following questions: 1. What are the structures of robust paths? 2. Given the structural information, can we find algorithms to trace the robust paths $\mathcal{P}(\mathcal{V})$ under different $\mathcal{V}$?

\textbf{Contribution.} 
In this paper, we answer both questions positively: 
\begin{itemize}
    \item[1.] We characterize the geometry of the robust paths as the Bregman projection of curves (whose geometries are defined by $\mathcal{V}$) onto the feasible region.
    \item[2.] Once the appropriate geometric lens is established, we find a surprisingly simple way to approximate robust paths. We connect the following two \emph{conceptually distinct}, yet \emph{geometrically similar} solution sets: 
    a) the robust paths of robust optimization problems (RC) and b) the optimization paths of the deterministic optimization counterparts (P). Specifically: 2.1. The {\ProximalPointMethodFull} ({\ProximalPointMethodAbb}) optimization paths for solving the deterministic counterpart (P) initialized at the ``most robust'' solution are close approximate, sometimes even exact, robust paths. 2.2. The design of the robust path uncertainty set shape, $\mathcal{V}$ is equivalent to the choice of the {\ProximalPointMethodAbb} distance-generating function; adjusting the cadence of the robust solutions' radii, $r$, corresponds exactly to adjusting the step-size of the {\ProximalPointMethodAbb}. 2.3. The distance between the {\ProximalPointMethodAbb} approximation and the exact robust path hinges precisely on the geometry of the feasible region $\mathcal{X}$ and the uncertainty set shape $\mathcal{V}$.
\item[3.] We computationally approximate robust paths in two settings---portfolio optimization and adversarial deep learning---to show the quality and efficiency of the resulting computational strategies.
\end{itemize}

\subsection{Related Works}
Calibrating the uncertainty set design $(r,\mathcal{V})$ to balance robustness and efficiency is a major challenge in applying robust optimization to real-world problems. \cite{Bertsimas2004} introduced the budgeted uncertainty set where the robustness and efficiency trade-off can be adjusted via a single perturbation budget hyperparameter. Under the assumption that the uncertainty follows a symmetric and bounded distribution, they provide probabilistic bounds on constraint violation as a theoretical remedy for setting the budget hyperparameter. A more recent line of work on data-driven distributionally robust optimization \citep{Gao_WDRO, mohajerin2018data, Blanchet_1} explicitly models the ambiguity in the uncertainty distribution and optimizes the worst-case objective over an ambiguity set of distributions within a fixed distance of the empirical distribution, \eg, a Wasserstein ball of radius $r$ centered at the empirical distribution. \cite{mohajerin2018data} provides theoretical guarantees that for a sufficiently large ambiguity set size, the ambiguity set contains the true distribution with high probability. However, the uncertainty/ambiguity set size and shape designed via probabilistic guarantees are empirically observed to be too conservative. Furthermore, in practice, many problems lack such statistical information to begin with for meaningful probabilistic guarantees. Thus, practitioners often resort to a computationally costly approach of computing multiple robust solutions under varying $(r,\mathcal{V})$, before selecting the robust solution with the best out-of-sample performance \citep{Ben_Phi_divergence, sim2021new, Robust_Satisficing, mohajerin2018data}. We quote \cite{Ben_Phi_divergence}: \emph{``However, it is a priori not possible to judge which uncertainty set is the `best'. We advocate the pragmatic approach to perform the robust optimization for different choices of uncertainty set, and then select the one that leads to the best optimal objective value.''}. Therefore, we argue that the challenge is in efficiently obtaining the entire robust path. In this work, we show that the robust path of (RC) can be approximated via a single {\ProximalPathFull} of (P), \emph{without} the need to solve the (RC) multiple times. We provide approximation error bounds that depend solely on the geometry of the uncertainty set and the feasible region.

Within the robust optimization literature, the properties of the robust paths have been studied in terms of their stability and continuity. \cite{Stability_Tim} and \cite{Stability_2} study the stability of the robust optimization optimal value and the robust optimal solution set with regard to variations of the uncertainty set. Our work has a different goal, which is the geometrical characterization of the \emph{entire} robust path and subsequently the algorithmic approximation of the robust paths of $(\mathrm{RC})$ via optimization paths of $(\mathrm{P})$. \cite{Iancu2014} studies another important solution set for robust optimization. They show that under a fixed uncertainty set design $(r,\mathcal{V})$, the set of robustly (worst-case) optimal solutions is not necessarily unique; subsequently, they introduce Pareto robustly optimal solutions within the set of robustly optimal solutions: if no other solutions perform at least as well for all uncertainty realizations and strictly better for some uncertainty realizations. Our robust path is a different robust solution set under varying $r$.

Robust optimization, under some conditions, has been shown to be equivalent to regularized optimization \citep{xu2009robustness,shafieezadeh2015distributionally, mohajerin2018data}. For this reason, the robust path is closely connected to the concept of regularization path in the statistics and machine learning literatures. A regularization path is defined as the set of optimal model weights under varying regularization hyperparameters. Existing works on regularization paths can be categorized as descriptive (analyzing the structural properties and the complexity of the regularization path) \citep{rp_d1,rp_d2,rp_d3} and constructive (algorithms generating an approximate or exact regularization path) \citep{rp_p0, rp_p1, rp_p2, rp_p3, rp_p4}. The connection between the first-order method trajectory and the regularization path (often termed explicit regularization of first-order methods) has been observed empirically \citep{barrett2021implicit} and studied theoretically \citep{opt_reg_path}. The regularization path literature focuses on \emph{unconstrained} regularized optimization problems encountered in learning. Our work differs from the literature and addresses \emph{constrained} robust/regularized problems. Most interestingly, under the constrained setting, we show that 1. the geometry of the robust/regularization path and 2. the quality of the approximation of robust/regularization paths of $(\mathrm{RC})$ by the optimization paths for $(\mathrm{P})$ both depend critically on the geometry of the feasible region $\mathcal{X}$. Our constrained results subsume the unconstrained setting in the case of linear objective functions.


\subsection{Structure of Paper}
\begin{figure}[ht]
\center
\begin{tikzpicture}[scale=0.8][ht]
\node at (1.25,1.4) {\small Robust Path};
\draw[rounded corners=2.5pt, thick] (0,0) rectangle (2.5,1) node[pos=.5] {$\mathcal{P}'(\mathcal{V})$};
\node at (10,1.4) {\small {\CentralPathFullCap}};
\draw[rounded corners=2.5pt, thick] (8.5,0) rectangle (11.5,1) node[pos=.5] { $\left\{ x_{\mathrm{{\CentralPathAbb}}}(\omega)\right\}$};
\node at (19.2,1.4) {\small {\ProximalPathFullCap}};
\draw[rounded corners=2.5pt, thick] (17.5,0) rectangle (20.5,1) node[pos=.5] { $\left\{ x_{k} \right\}$};

\node at (5.5,2.75) {\small \fbox{
$\begin{array}{c}
 \text{Approximation error bound} \\ 
 \text{between} \\  
 \text{$\mathcal{P}'(\mathcal{V})$ and $\left\{ x_{\mathrm{{\CentralPathAbb}}}(\omega)\right\}$}   \\ 
 \end{array}$}};
\draw[<->,decorate,
        decoration={snake, pre length=10pt, post length=10pt},thick] (2.5,0.8) -- (8.5,0.8) node[midway,above] {\small \textbf{Theorem \ref{theorem: uniform bound CP and RP}}};

\node at (5.5,-1.5) {\small \fbox{$\begin{array}{c}
\text{Sufficient condition for} \\  \mathcal{P}'(\mathcal{V})=\{x_\mathrm{{\CentralPathAbb}}(\omega)\} \\
\end{array}$}}; 
\draw[<->,thick] (2.5,0.2) -- (8.5,0.2) node[midway,below] {\small\textbf{Corollary \ref{corollary: BPP is RP}} };

\node at (14.5,2.75) {\small \fbox{
$\begin{array}{c}
 \text{Approximation error bound} \\ 
 \text{between} \\  
 \text{$\left\{ x_{\mathrm{{\CentralPathAbb}}}(\omega)\right\}$ and $\{x_k\}$}   \\ 
 \end{array}$}};
\draw[<->,decorate,
        decoration={snake, pre length=10pt, post length=10pt},thick] (11.5,0.8) -- (17.5,0.8) node[midway,above] {\small \textbf{Theorem \ref{theorem: PPM as approx BPP}}};

\node at (14.5,-1.75) {\small \fbox{$\begin{array}{c}
\text{$ 
\{x_{k}\} \text{ is monotone on } \mathcal{X} \ \Rightarrow \  \{x_{k}\} = \{x_{\mathrm{{\CentralPathAbb}}}\left(\omega_{k}\right)\}$} \\
\\
\text{$\text{$\mathcal{X}$ and $\mathcal{V}$ are polar pairs} \ \Rightarrow \  \{x_{k}\} \subset  \{x_{\mathrm{{\CentralPathAbb}}}\left(\omega_{k}\right)\}$}
\end{array}$}}; 
\draw[<->,thick] (11.5,0.2) -- (17.5,0.2) node[midway,below,align=center] {\small \textbf{Proposition \ref{prop: monotone_BPP_is_PPM}},  \textbf{Proposition \ref{prop: X_matches_V}}};
\draw[<->,thick] (1,0) -- (1,-3.5) -- (19.5,-3.5) -- (19.5,0);
\node at (10.25,-3.95) {\small \textbf{Theorem \ref{theorem: PPMP_Exact_RP}}: Sufficient condition for $\mathcal{P}'(\mathcal{V})=\{x_k\}$ };



\end{tikzpicture}
\caption{Summary of theory in this paper (Sections \ref{section: RP_geometry} and \ref{section: approxiamte_robust_path}). 
Theorem \ref{theorem: RP_CP_PP are Bregman Projections} reveals a geometric lens that enables this whole analysis. Wavy lines indicate approximation bounds; solid lines indicate equivalence results. $\mathcal{P}'(\mathcal{V})$ (Definition \ref{definition: unique robust path}) is a subset of $\mathcal{P}(\mathcal{V})$ (Definition \ref{def:robust_path}) that contains a unique solution under each uncertainty radius $r$.}
\label{fig: summary_sec_2_results}
\end{figure}
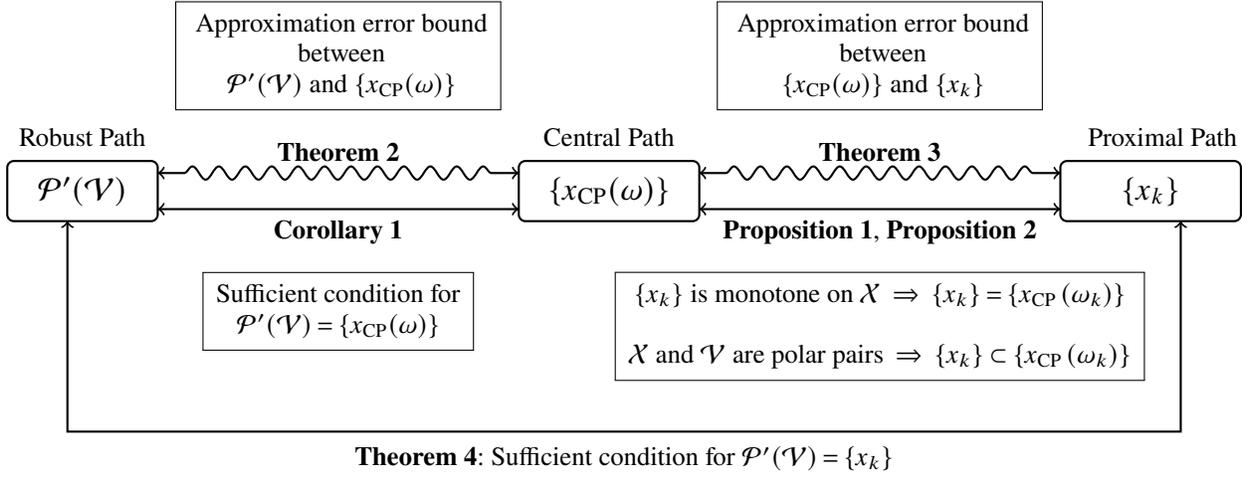

The presentation of our results is structured as follows. Section \ref{section: Preliminary} introduces key technical concepts fundamental to our analysis. In Section \ref{section: RP_geometry}, we present a geometric lens that unifies the characterization of the robust paths of $(\mathrm{RC})$ and two optimization paths for $(\mathrm{P})$, all as Bregman projections of curves onto the feasible region.  In Section \ref{section: approxiamte_robust_path}, we show the robust paths of $(\mathrm{RC})$ can be algorithmically approximated via the optimization paths for solving $(\mathrm{P})$. The structure of the main results is summarized in Figure \ref{fig: summary_sec_2_results}. We derive a sharp approximation error bound that depends on the geometry of the feasible region and the uncertainty set. We also present two special cases where the approximation error is zero based on polyhedral analysis and duality. In Section \ref{section: application}, we validate our theory in two settings: portfolio optimization and adversarial deep learning.


\section{Preliminary}
\label{section: Preliminary}
In this section, we review concepts that are fundamental to our analysis in Sections \ref{section: RP_geometry} and \ref{section: approxiamte_robust_path}.
\subsection{Notations}
 We adopt the notations of \cite{Rockafellar_Convex_Analysis}. We denote the extended real numbers as $\overline{\mathbb{R}}=\mathbb{R}\cup \{-\infty, +\infty\}$, the nonnegative real numbers as $\mathbb{R}_{+}=[0,+\infty]$. Given an extended real-valued function $f: \mathbb{R}^{n} \rightarrow \overline{\mathbb{R}}$, its \textit{effective domain} is defined as $\mathrm{dom}(f)=\{x\in \mathbb{R}^{n}: f(x) < +\infty\}$, its \textit{epigraph} is defined as $\mathrm{epi}(f)=\{(x,t)\in \mathbb{R}^{n} \times \mathbb{R}: f(x)\leq t\}$. $f$ is \textit{proper} if $f(x) > -\infty$ for all $x\in \mathbb{R}^{n}$ and there exists $x\in \mathbb{R}^{n}$ such that $f(x)< +\infty$, \ie, $\mathrm{dom}(f) \neq \varnothing$. In addition, $f$ is \textit{closed} if $\mathrm{epi}(f)$ is closed. Let $\mathcal{C}\subseteq \mathbb{R}^{n}$ be a nonempty closed convex set, we denote the \textit{interior} of $\mathcal{C}$ as $\mathrm{int}(\mathcal{C})$, the \textit{relative interior} of $\mathcal{C}$ as $\mathrm{ri}(\mathcal{C})$, the \textit{boundary} of $\mathcal{C}$ as $\mathrm{bd}(\mathcal{C})$ and the \textit{affine hull} of $\mathcal{C}$ as $\mathrm{Aff}(\mathcal{C})$. The \textit{polar set} of $\mathcal{C}$ is defined as $\mathcal{C}^{\circ}=\{y\in \mathbb{R}^{n}: \langle x,y\rangle\leq 1, \forall x \in \mathcal{C}\}$. In addition, we define the \textit{support function} of $\mathcal{C}$ as $\sigma(\ \cdot \ \vert  \mathcal{C}) = \sup_{y\in \mathcal{C}} \langle \ \cdot \ , y \rangle$. Let $\mathcal{V}$ be a closed and convex set containing the origin, we define the \textit{gauge function} induced by $\mathcal{V}$ as $\|v\|_{\mathcal{V}}=\inf \{t \geq 0: v\in t\mathcal{V} \}$. The gauge function has the well-known dual formulation $\|v\|_{\mathcal{V}}=\sigma(v \vert \mathcal{V}^{\circ})$ \cite[Theorem~14.5]{Rockafellar_Convex_Analysis}. 

\subsection{Smoothness and Strict Convexity of Convex Sets}
\begin{definition}[Strict convexity and Smoothness of Convex Sets]
    A convex set $\mathcal{C}$ with $\mathrm{int}(\mathcal{C}) \neq \varnothing$ is \textit{strictly convex} if for any distinct $x$ and $y \in \mathrm{bd}(\mathcal{C})$, $\{\alpha x+ (1-\alpha)y : \alpha \in (0,1)\} \subset \mathrm{int}(\mathcal{C})$. A point $x\in \mathrm{bd}(\mathcal{C})$ is \textit{regular} if the supporting hyperplane of $\mathcal{C}$ at $x$ is unique, $\mathcal{C}$ is \textit{smooth} if $x$ is regular for all $ x \in \mathrm{bd}(\mathcal{C})$.
\end{definition}

\subsection{Legendre Functions}
We mostly work with Legendre functions as defined in \cite{Rockafellar_Convex_Analysis}, Section 26. 
\begin{definition}[Legendre]
    Suppose $\psi: \mathbb{R}^{n} \rightarrow \mathbb{R} \cup \{+\infty\}$ is proper, closed and convex. The function $\psi$ is \textit{Legendre} (or a \textit{Legendre function}), if $\psi$ is both essentially smooth and essentially strictly convex, \ie, $\psi$ satisfies
    \begin{itemize}
        \item[i)] $\mathrm{int}(\mathrm{dom}(\psi)) \neq \varnothing$;
        \item[ii)] $\psi$ is differentiable on $\mathrm{int}(\mathrm{dom}(\psi))$;
        \item[iii)] $\|\nabla \psi(x_{n})\|_{2} \rightarrow + \infty $, for any sequence $(x_{n}) \subset \mathrm{int}(\mathrm{dom}(\psi))$ such that $x_{n} \rightarrow x \in \mathrm{bd}(\mathrm{dom}(\psi))$;
        \item[iv)] $\psi$ is strictly convex on $\mathrm{int}(\mathrm{dom}(\psi))$.
    \end{itemize}
\end{definition}
In particular, the following property of Legendre functions is useful in studying the robust path and optimization paths in both the primal and dual space as defined by a bijection $\nabla \psi$:
\begin{lemma}[\cite{Rockafellar_Convex_Analysis}, Theorem 26.5]
\label{lemma: legendre_bijective}
Let a proper closed convex function $\psi: \mathbb{R}^{n} \rightarrow \mathbb{R} \cup \{+\infty\}$ be a Legendre function. Then $\nabla \psi: \mathrm{int}(\mathrm{dom}(\psi)) \rightarrow \mathrm{int}(\mathrm{dom}(\psi^{*}))$ is a bijection, with $(\nabla \psi)^{-1} = \nabla \psi^{*}:\mathrm{int}(\mathrm{dom}(\psi^{*})) \rightarrow \mathrm{int}(\mathrm{dom}(\psi)) $.
\end{lemma}

\subsection{Bregman Divergence, Bregman Projection and Generalized Proximal Operator}
Bregman divergence provides a generalized notion of proximity between two points. 
\begin{definition}[Bregman Divergence]
\label{def: Bregman Divergence}
    Let the distance-generating function (d.g.f.), $\psi: \mathbb{R}^{n} \rightarrow \mathbb{R} \cup \{+\infty\}$ be proper closed convex and differentiable on $\mathrm{int}(\mathrm{dom}(\psi))$, the \textit{Bregman divergence} induced by $\psi$, $D_{\psi}(\cdot ,\cdot): \mathbb{R}^{n} \times \mathrm{int}(\mathrm{dom}(\psi)) \rightarrow \mathbb{R}_{+}$ is defined as: 
    $$ D_{\psi}(x,y) := \psi(x) - \psi(y) - \langle \nabla \psi(y), x-y \rangle. $$
\end{definition}

Bregman projection generalizes the Euclidean projection. 
\begin{definition}[Bregman Projection]
Fix a function $\psi: \mathbb{R}^{n} \rightarrow \mathbb{R} \cup \{+\infty\}$ that is proper closed  convex and differentiable on $\mathrm{int}(\mathrm{dom}(\psi))$, and a closed convex set $\mathcal{S}\subset \mathrm{int}(\mathrm{dom}(\psi))$. We define the \textit{Bregman projection} associated with $\psi$ of a point $y \in \mathrm{int}(\mathrm{dom}(\psi))$ onto $\mathcal{S}$ as: $\Pi^{\psi}_{\mathcal{S}}(y) = \argmin_{x\in \mathcal{S}} D_{\psi} \left(x,y\right).$
\end{definition}

A Legendre $\psi$ ensures the existence and the uniqueness of the Bregman projection.
\begin{lemma}[\cite{bauschke1997legendre}, Theorem 3.12]
If $\psi$ is Legendre, then $\Pi^{\psi}_{\mathcal{S}}(y)$ is a singleton.
\end{lemma}

Bregman projection has the following variational characterization.
\begin{lemma}[\cite{bauschke1997legendre}, Proposition 3.16]\label{lemma: Bregman Projection VI} 
A point $x'\in \mathcal{C}$ is $\Pi^{\psi}_{\mathcal{C}}(x)$ if and only if
\begin{equation}
\label{eq: proj_lemma}
\langle \nabla \psi(x) - \nabla \psi(x'), y-x'  \rangle \leq 0 , \quad \forall y \in \mathcal{C}.
\end{equation}
\end{lemma}

The Bregman projection satisfies the following property.
\begin{lemma}
\label{lemma: Breg_Proj_calc_1}
    Let $\psi$ be a Legendre function, let $\mathcal{C}$ be a closed convex set, then the associated Bregman projection operator satisfies
    $$\Pi^{\psi}_{\mathcal{C}} = \Pi^{\psi}_{\mathcal{C}} \circ \Pi^{\psi}_{\mathrm{Aff}(\mathcal{C})}$$
\end{lemma}
\begin{proof}{Proof.}
 We defer the proof to Appendix \ref{APP: proof of Breg_Proj_calc_1}.
\end{proof}

The Bregman projection onto a closed convex cone has the following dual characterization.
\begin{lemma}[\cite{bauschke2003duality}, Theorem 3.1]
\label{lemma:dual_cone bregman_proj}
    Let $\psi$ be a Legendre function, let $\mathcal{K}$ be a closed convex cone in $\mathrm{dom}(\psi)$ with its dual cone denoted as $\mathcal{K}^{*}$, then
    $$ \nabla \psi \left( \Pi^{\psi}_{\mathcal{K}+x_{0}}\left(y_{0}\right) \right) = \Pi^{\psi^{*}}_{\mathcal{K}^{*}+\nabla \psi(y_{0})}\left(\nabla \psi(x_0) \right). $$
\end{lemma}

The Bregman proximal operator generalizes the standard proximal operator by replacing its Euclidean distance with Bregman divergence.
\begin{definition}[Bregman Proximal Operator]
\label{def: proximal operator}
Fix a function $\psi: \mathbb{R}^{n} \rightarrow \mathbb{R} \cup \{+\infty\}$ that is proper closed  convex and differentiable on $\mathrm{int}(\mathrm{dom}(\psi))$. Let $\eta > 0$. We define the \textit{Bregman proximal operator} of a linear function $f(\cdot) = \langle c, \cdot \rangle$ and a closed convex set $\mathcal{S}\subset \mathrm{int}(\mathrm{dom}(\psi))$ induced by $\psi$ as: $ \mathrm{Prox}^{\psi}_{c,\mathcal{S}}(y,\eta) = \argmin_{x\in \mathcal{S}} \langle c,x \rangle + \frac{1}{\eta} D_{\psi}(x, y). $
\end{definition}


\section{Unified Geometric View of Robust Path and Two Optimization Paths}
\label{section: RP_geometry}

In this section, we develop a unified Bregman projection view that allows us to geometrically relate the robust path of the robust counterpart $(\mathrm{RC})$, and what we will define as optimization paths for solving the deterministic counterpart $(\mathrm{P})$.

\subsection{Definitions of Two Optimization Paths} 
We introduce two optimization paths for the deterministic optimization problem $(\mathrm{P})$ as follows. In Definition \ref{def:proximal-path}, we utilize the proximal point method, which is a fundamental tool in the analysis of modern optimization algorithms. 
\begin{definition}[Bregman Proximal Point Method and Bregman Proximal Point Path]\label{def:proximal-path}
Given a distance-generating function $\psi$, a step-size sequence $\{\lambda_{k}> 0 \}$, and a starting point $x_{0}\in \mathcal{X}$, the \emph{Bregman proximal point method} for solving problem $(\mathrm{P})$ generates a sequence $\{x_{k}\}$ 
$$x_{k+1} = \argmin_{x\in \mathcal{X}} \langle a_0,x  \rangle + \lambda_{k}D_{\psi}(x,x_{k}), \, k = 0, 1, \dots$$
The \emph{Bregman proximal point path} is defined as the sequence $\{x_{k}\}$.
\end{definition}
Intuitively, $x_{k+1}$ is a point in $\mathcal{X}$ that tries to minimize $\langle a_0, x \rangle$ without deviating too much from anchor point $x_k$, \ie, the previous iterate. The next definition describes points that try to minimize $\langle a_0, x \rangle$ without deviating too much from some \emph{fixed} anchor point $x_0$. 
\begin{definition}[Bregman Central Path]\label{def:central-path}
The \emph{Bregman central path} for solving problem  $(\mathrm{P})$, induced by distance-generating function $\psi$ and initialized at $x_{0}$ is the set $\{x_{\mathrm{{\CentralPathAbb}}}(\omega): \omega \in [0,\infty)\}$ defined as
$$x_{\mathrm{{\CentralPathAbb}}}(\omega) = \argmin_{x\in \mathcal{X}} \langle a_0,x  \rangle + \omega D_{\psi}(x,x_{0}). $$

Throughout this paper, we work with strictly convex $D_{\psi}$ (Assumption \ref{assumption: blanket varphi} below). Therefore the optimal solutions are unique and we use ``$=$'' instead of ``$\in$'' in Definitions \ref{def:proximal-path} and \ref{def:central-path}.

For ease of notation, we henceforth refer to the Bregman proximal point path as the \emph{{\ProximalPathFull}} ({{\ProximalPathAbb}}); and the Bregman central path as the \emph{{\CentralPathFull}} ({\CentralPathAbb}). We abuse the notion of central path in this paper: we assume $D_{\psi}(\cdot, x_0)$ is not necessarily a self-concordant barrier function of the feasible region $\mathcal{X}$, hence the central path is not necessarily strictly contained in the interior of the feasible region $\mathcal{X}$.

\end{definition}

\subsection{The Geometric View}
We make the following assumptions throughout the remainder of the paper.
\begin{assumption}
\label{assumption: blanket V}
The uncertainty set $\mathcal{V}$ is compact, smooth, and strictly convex, with $0 \in \mathrm{int}(\mathcal{V})$.
\end{assumption}

\begin{assumption}
\label{assumption: blanket varphi}
Define $\varphi(\cdot)=g \circ \|\cdot\|_{\mathcal{V}^{\circ}}$, where $g:\mathbb{R} \rightarrow \mathbb{R}_{+}$ is Legendre, $g(0)=0$ and $\nabla g(0) = 0$.
\end{assumption}
Note that while Assumption \ref{assumption: blanket V} places a requirement on the uncertainty set of $(\mathrm{RC})$, Assumption \ref{assumption: blanket varphi} is a technical assumption to facilitate proofs, not a restriction on $(\mathrm{RC})$.

In Definition \ref{def:robust_path}, we defined robust path as the set of (potentially non-unique) optimal solutions of $(\mathrm{RC})$ with different radii. Now with the help of Assumptions \ref{assumption: blanket V} and \ref{assumption: blanket varphi}, we refine Definition \ref{assumption: blanket V} into Definition \ref{definition: unique robust path}, where each radius $r$ corresponds to a \emph{unique} optimal solution. The main implication of the solution uniqueness of $\mathcal{P}'(\mathcal{V})$ is that its geometry can be precisely characterized via Bregman projection.
\begin{definition}[(Characterizable) Robust Path]
\label{definition: unique robust path}
The \emph{(characterizable) robust path} of $(\mathrm{RC})$ under $\mathcal{V}$ is defined as
    $ \mathcal{P}'(\mathcal{V}) = \left\{x_{\mathrm{R}}'(\omega, \mathcal{V}):   \omega\in [0,\infty) \right\},$
    where $$x_{\mathrm{R}}'(\omega) := x_{\mathrm{R}}'(\omega, \mathcal{V}) = \argmin_{x\in \mathcal{X}} \langle a_{0}, x \rangle + \omega \varphi(x) = \argmin_{x\in \mathcal{X}} \langle a_{0}, x \rangle  + \omega \cdot g \circ \|x\|_{\mathcal{V}^{\circ}}.$$
The robust solution $x_{\mathrm{R}}$ and the efficient solution $x_{\mathrm{E}}$ are defined as
$$
    x_{\mathrm{R}} = \lim_{\omega \rightarrow \infty} x_{\mathrm{R}}'(\omega, \mathcal{V}) = \argmin_{x\in \mathcal{X}} \varphi(x) \quad \text{and} \quad x_{\mathrm{E}} = \lim_{\omega \rightarrow 0} x_{\mathrm{R}}'(\omega, \mathcal{V})  = \argmin_{x\in \mathcal{X}} \langle a_{0},x \rangle.
$$
\end{definition}
Note that $x_{\mathrm{R}}'(\omega, \mathcal{V})$ is defined as the solution of a regularized problem, not directly of a robust optimization problem. This is possible due to the following duality results. 

\begin{lemma}
\label{lemma: dual_RC}
    The dual problem of $(\mathrm{RC})$ under $\Xi(r,\mathcal{V}) = \left\{\xi \in \mathbb{R}^{n}:    \|\xi\|_{\mathcal{V}} \leq r  \right\}$ is 
    \begin{equation}
    \label{eq: dual_RC}
    (\mathrm{RCD}) \quad \min_{x\in \mathcal{X}} \langle a_{0}, x \rangle + r\|x\|_{\mathcal{V}^{\circ}}.
    \end{equation}  
In addition, strong duality holds under $0 \in \mathrm{int}(\mathcal{V})$.
\end{lemma}

\begin{lemma}
\label{lemma: unique robust path}
The following statements are true:
\begin{itemize}
    \item[(i)] The path in Definition \ref{definition: unique robust path} is indeed a set of robust solutions: $\mathcal{P}'(\mathcal{V}) \subset \mathcal{P}(\mathcal{V})$.
    \item[(ii)] Given any solution $x_{\mathrm{R}}'(\omega, \mathcal{V})$ on the path corresponding to a specific regularization strength $\omega$, the corresponding robust optimization uncertainty set radius can be identified. In other words, $$x_{\mathrm{R}}'(\omega, \mathcal{V}) \in \argmin_{x \in \mathcal{X}}  \max_{\xi  \in\Xi(r(\omega),\mathcal{V})} \langle a_{0} +  \xi,x \rangle $$ where $r(\omega) = \omega  \nabla g (\|x_{\mathrm{R}}'\left(\omega, \mathcal{V})\|_{\mathcal{V}^{\circ}} \right)$.
\end{itemize}
\end{lemma}
We defer the proofs of Lemmas \ref{lemma: dual_RC} and \ref{lemma: unique robust path} to Appendix \ref{sec: proof proposition: unique robust path}.

Now we are ready for the first main result. Theorem \ref{theorem: RP_CP_PP are Bregman Projections} reveals that the robust path $\mathcal{P}'(\mathcal{V})$ of $(\mathrm{RC})$ and two appropriately defined optimization paths of $(\mathrm{P})$ are in fact geometrically similar under the lens of Bregman projection.
\begin{theorem}
\label{theorem: RP_CP_PP are Bregman Projections}
Denote $\mathcal{P}'(\mathcal{V})$ as the \emph{robust path} of $(\mathrm{RC})$ according to Definition \ref{definition: unique robust path}.
Define $x_{\mathrm{R}} = \lim_{\omega \rightarrow \infty} x_{\mathrm{R}}'(\omega,\mathcal{V}) = \argmin_{x\in \mathcal{X}} \varphi(x)$. 
Denote $\left\{x_{\mathrm{{\CentralPathAbb}}}(\omega) : \omega \in [0,\infty) \right\}$ as the \emph{{\CentralPathFull}} for solving $(\mathrm{P})$, using $\varphi(\cdot)=g \circ \|\cdot\|_{\mathcal{V}^{\circ}}$ as the d.g.f. and initialized at $x_{\mathrm{R}}$. 
Denote $\left\{ x_{k} \right\}$ as the \emph{{\ProximalPathFull}} for solving $(\mathrm{P})$, using $\varphi(\cdot)=g \circ \|\cdot\|_{\mathcal{V}^{\circ}}$ as the d.g.f., under a step-size sequence $\{\lambda_{k}> 0 : \, \sum_{k=0}^{\infty} \lambda_{k}^{-1} = \infty \}$ and initialized at $x_{\mathrm{R}}$.
 The three paths have the following Bregman projection interpretation:
\\
\\
$
\begin{aligned}
  \text{\emph{(Robust Path)}}\hspace{3.0cm}  \mathcal{P}'(\mathcal{V})&= \left\{  \Pi^{\varphi}_{\mathcal{X}} \left( \nabla \varphi^{*}\left(\nabla\varphi( \ 0 \ ) - \omega^{-1} a_{0} \right) \right) : \omega \in [0,\infty) \right\},\\
  \text{\emph{({\CentralPathFullCap})}}\quad \left\{x_{\mathrm{{\CentralPathAbb}}}(\omega) : \omega \in [0,\infty) \right\}\!&= \left\{\Pi^{\varphi}_{\mathcal{X}} \left( \nabla \varphi^{*}\left(\nabla\varphi(x_{\mathrm{R}}) - \omega^{-1} a_{0} \right) \right) : \omega \in [0,\infty) \right\},\\
  \text{\emph{({\ProximalPathFullCap})}}\hspace{3.16cm} x_{k+1} \!&= \; \; \, \Pi^{\varphi}_{\mathcal{X}} \left( \nabla \varphi^{*}\left(\nabla\varphi(x_{k}) - \lambda_{k}^{-1} \ a_{0} \right) \right) , \; \; \text{for } k  = 0,1,\cdots, \ x_{0} = x_{\mathrm{R}}.
\end{aligned}
$
\\
\end{theorem} 

\textit{Remark.} 
Theorem \ref{theorem: RP_CP_PP are Bregman Projections} reveals that the robust path, the central path, and the proximal path can be characterized as the Bregman projection of a curve (induced by a ray in the dual space $\nabla \varphi$) onto the feasible set $\mathcal{X}$. More specifically, the robust path and the {\CentralPathFull} only differ in their respective initial points of the curves, \ie, $\nabla \varphi(0)$ and $\nabla \varphi(x_{\mathrm{R}})$. Thus, as one may expect and as we will show in Section \ref{section: approxiamte_robust_path}, {\CentralPathFull}s are approximate, sometimes exact robust paths. Further, the ``distance'' between the {\CentralPathFull} and the robust path depends on the ``distance'' between the robust solution $x_{\mathrm{R}}$ and the origin $0$ in a precise and non-trivial manner. We also point out that the {\CentralPathFull} can be viewed as the Bregman projection of a curve onto $\mathcal{X}$, while the {\ProximalPathFull} is generated via successive Bregman projection of small curve segments onto $\mathcal{X}$. In Section \ref{section: approxiamte_robust_path}, we show the {\ProximalPathFull}s are approximate, sometimes exact {\CentralPathFull}s. Finally, leveraging the central paths as intermediaries, we show that the proximal paths are approximate, under some conditions, exact robust paths.

\begin{figure}[ht]
    \centering
    \begin{subfigure}[t]{0.48\linewidth}
        \centering
        \includegraphics[width=\linewidth]{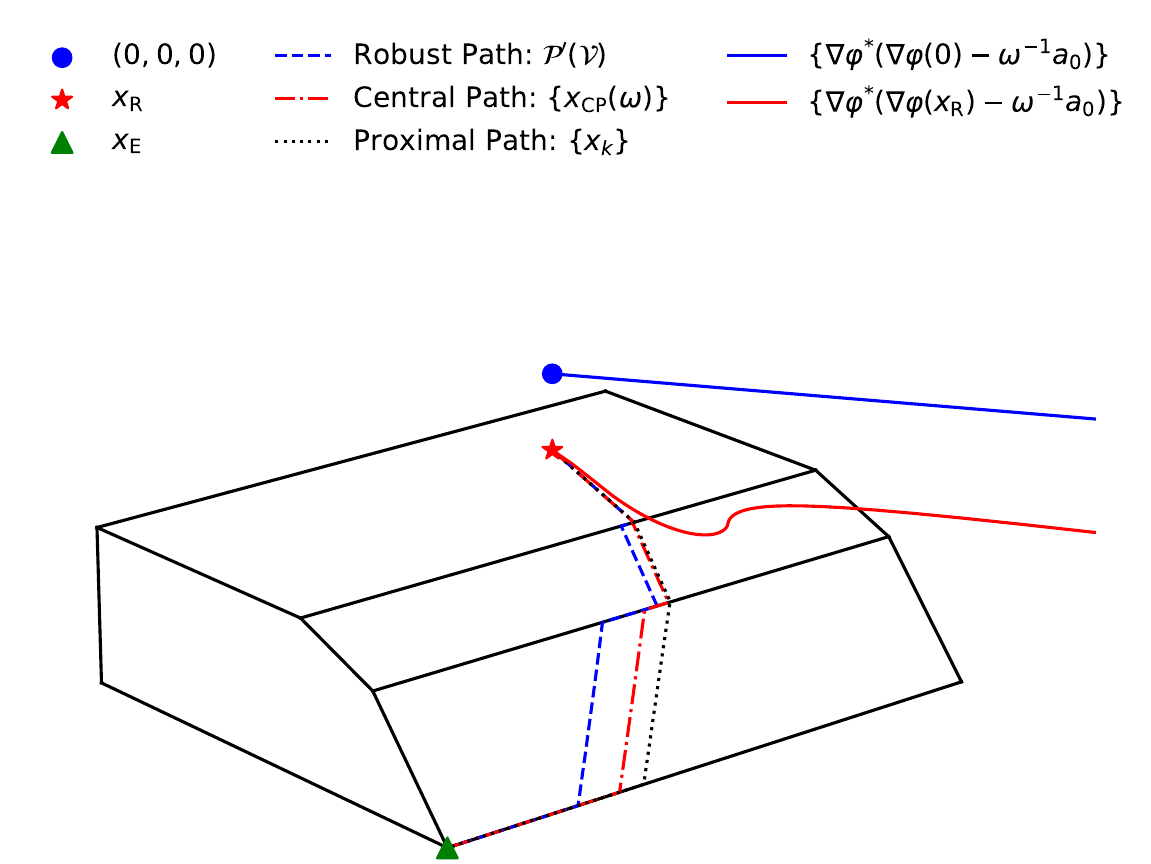}
        \caption{\footnotesize Geometry of the robust path $\mathcal{P}'(\mathcal{V})$, the {\CentralPathFull} $\{x_{\mathrm{{\CentralPathAbb}}}(\omega)\}$, and the {\ProximalPathFull} $\{x_{k}\}$. $\mathcal{V} = \{x\in \mathbb{R}^{3}: \|x\|_{p=5/3} \leq 1\}$, $\varphi = \|x\|_{\mathcal{V}^{\circ}}^{2} = \|x\|_{q=5/2}^{2}$.}
        \label{fig: unified_geometric_view_a}
    \end{subfigure} \hfill
    \begin{subfigure}[t]{0.46\linewidth}
        \centering
        \includegraphics[width=\linewidth]{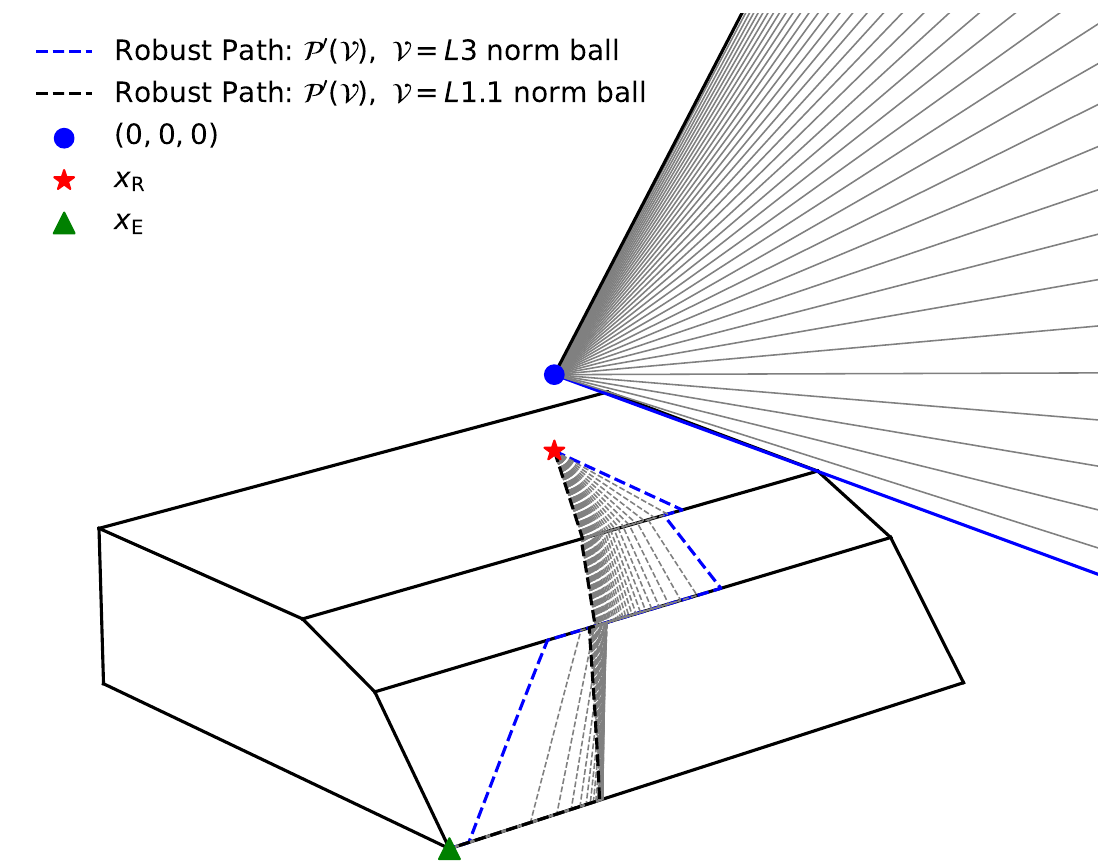}
        \caption{\footnotesize  Geometry of a set of robust paths under varying uncertainty sets. $\{\mathcal{P}'(\mathcal{V}): \: \mathcal{V} = \{x\in \mathbb{R}^{3}: \|x\|_{p} \leq 1\}, \ p \in [1.1, 3]\}$, $\varphi = \|x\|_{\mathcal{V}^{\circ}}^{2}.$}
        \label{fig: unified_geometric_view_b}
    \end{subfigure}
    \caption{Theorem \ref{theorem: RP_CP_PP are Bregman Projections}: Unified Bregman projection view of robust paths of $(\mathrm{RC})$, {\CentralPathFull} of $(\mathrm{P})$, and {\ProximalPathFull} of $(\mathrm{P})$.}
    \label{fig: unified_geometric_view}
\end{figure}

In Figure \ref{fig: unified_geometric_view_a}, we visualize, as Theorem \ref{theorem: RP_CP_PP are Bregman Projections} entails, the Bregman projection interpretation of robust path $\mathcal{P}'(\mathcal{V})$ of $\mathrm{(RC)}$ together with the two optimization paths: {\CentralPathFull} $\{x_{\mathrm{{\CentralPathAbb}}}(\omega)\}$ and {\ProximalPathFull} $\{x_{k}\}$ of $(\mathrm{P})$ for a single $\mathcal{V}$. The distance-generating function $\varphi$ is induced by $\mathcal{V}$ (\ie, $\varphi(\cdot)=g \circ \|\cdot\|_{\mathcal{V}^{\circ}}$). Figure \ref{fig: unified_geometric_view_b} presents the geometry of a set of robust paths $\mathcal{P}'(\mathcal{V})$ under varying $\mathcal{V}$ designs.

\subsection{Proof of Theorem \ref{theorem: RP_CP_PP are Bregman Projections}}

We prove the representation of robust path in Section \ref{subsubsection-Thm1-proof-robust-path}, and prove the representations of central path and proximal path in Section \ref{subsubsection-Thm1-proof-CP-PP}.

\subsubsection{Robust Path: Bregman Projection and Dual Space.}\label{subsubsection-Thm1-proof-robust-path}

We prove the robust path's representation in Theorem \ref{theorem: RP_CP_PP are Bregman Projections} via proving three claims. In Claim \ref{claim: gauge_set_to_Legendre_dgf}, we show that when the uncertainty set and the post-composition function $g$ are nice, they induce a nice (Legendre) distance-generating function $\varphi$. In Claim \ref{claim: prox_is_proj}, we show that with a Legendre d.g.f., the Bregman proximal operator is a Bregman projection. In Claim \ref{claim: robust path is proj}, we show that since the robust path as defined in Definition \ref{definition: unique robust path} satisfies these properties, it can be represented as a Bregman projection of a curve induced by a ray in the dual space $\nabla \varphi$.

\begin{claim}
\label{claim: gauge_set_to_Legendre_dgf}
If $\mathcal{V}$ and $\varphi(\cdot)=g \circ \|\cdot\|_{\mathcal{V}^{\circ}}$ satisfy Assumption \ref{assumption: blanket V} and Assumption \ref{assumption: blanket varphi} respectively, we have
\begin{itemize}
    \item[(A) Well-defined:]  $\mathcal{X} \subseteq \mathrm{dom}(\varphi)$, \ie, $\varphi$ is well-defined over $\mathcal{X}$. 
    \item[(B) Legendre:] $\varphi$ is Legendre.
    \item[(C) Zero at Origin:] $\varphi(0)=0$ and $\nabla \varphi (0) = 0$.
\end{itemize}
\end{claim}

\begin{proof}{Proof.}
\emph{(A) Well-defined.} $\mathrm{dom}(\varphi) = \mathrm{dom}(\|\cdot\|_{\mathcal{V}^{\circ}}) = \Pi_{\mathbb{R}^{n}}(\mathrm{epi}(\|\cdot\|_{\mathcal{V}^{\circ}})) = \Pi_{\mathbb{R}^{n}}(\mathrm{cone}(\{x\in \mathbb{R}^{n}:\|x\|_{\mathcal{V}^{\circ}}\leq 1 \} \times \{1\})) = \mathbb{R}^{n} \supseteq \mathcal{X}$, where $\mathrm{cone}(\cdot)$ denotes the conic hull of a set. The second and third equalities are classic results in convex analysis, the last equality is by $0\in \mathrm{int}(\mathcal{V}^{\circ})$ (by \cite{Rockafellar_Convex_Analysis}, Corollary 14.5.1. and the assumption that $\mathcal{V}$ is compact and convex).

\emph{(B) Legendre: Essentially strictly convex.} Consider two points $x_{1}$ and $x_{2} \in \mathrm{int}(\mathrm{dom}(\varphi)) = \mathbb{R}^{n}$, denote $z = \beta x_{1} + (1-\beta) x_{2}$, for some $\beta\in (0,1)$. We show $\varphi(z) < \beta \varphi(x_{1}) + (1-\beta) \varphi(x_{2})$ for two cases: $\|x_{1}\|_{\mathcal{V}^{\circ}} = \|x_{2}\|_{\mathcal{V}^{\circ}}$ and  $\|x_{1}\|_{\mathcal{V}^{\circ}} \neq \|x_{2}\|_{\mathcal{V}^{\circ}}$.

Case One $\|x_{1}\|_{\mathcal{V}^{\circ}} =\|x_{2}\|_{\mathcal{V}^{\circ}}=\alpha$: Given $\mathcal{V}$ is closed, convex and smooth with $0 \in \mathcal{V}$, $\mathcal{V}^{\circ}$ is closed and strictly convex \cite[Proposition 3.2.7]{hiriart2004fundamentals}. Consequently, the sublevel set $\mathcal{L}_{\alpha} = \{x \in \mathbb{R}^{n}: \|x\|_{\mathcal{V}^{\circ}} \leq \alpha  \} = \alpha \cdot \mathcal{V}^{\circ}$ is strictly convex, hence $z \in \mathrm{int}(\mathcal{L}_{\alpha})$. In addition, given $\|\cdot\|_{\mathcal{V}^{\circ}}$ is continuous and homogeneous of degree of 1, we have $\|z\|_{\mathcal{V}^{\circ}} < \alpha$ or equivalently $\|z\|_{\mathcal{V}^{\circ}} < \beta \|x_{1}\|_{\mathcal{V}^{\circ}} + (1-\beta)\|x_{2}\|_{\mathcal{V}^{\circ}}$. Using the fact that $g$ is essentially strictly convex with $\nabla g(0)=0$, $g$ is monotonically increasing over $\mathbb{R}_{+}$. Therefore we have: $g\circ\|z\|_{\mathcal{V}^{\circ}} <  g(\beta \|x_{1}\|_{\mathcal{V}^{\circ}} + (1-\beta)\|x_{2}\|_{\mathcal{V}^{\circ}}) =\beta g \circ \|x_{1}\|_{\mathcal{V}^{\circ}} + (1-\beta)g\circ\|x_{2}\|_{\mathcal{V}^{\circ}}.$
    
Case Two $\|x_{1}\|_{\mathcal{V}^{\circ}} \neq \|x_{2}\|_{\mathcal{V}^{\circ}}$: Without loss of generality, we assume $\|x_{1}\|_{\mathcal{V}^{\circ}} < \|x_{2}\|_{\mathcal{V}^{\circ}}$. By the convexity of $\|\cdot \|_{\mathcal{V}^{\circ}}$: $\|z\|_{\mathcal{V}^{\circ}} \leq \beta\|x_{1}\|_{\mathcal{V}^{\circ}} + (1-\beta)\|x_{2}\|_{\mathcal{V}^{\circ}}.$ 
Given $g$ is monotonically increasing over $\mathbb{R}_{+}$ and essentially strictly convex: $g\circ \|z\|_{\mathcal{V}^{\circ}} \leq  g(\beta\|x_{1}\|_{\mathcal{V}^{\circ}} + (1-\beta)\|x_{2}\|_{\mathcal{V}^{\circ}})<  \beta g \circ \|x_{1}\|_{\mathcal{V}^{\circ}} + (1-\beta) g\circ \|x_{2}\|_{\mathcal{V}^{\circ}}. $

\emph{(B) Legendre: Differentiable.} Towards proving $\varphi$ is essentially smooth, we begin by showing $\varphi$ is differentiable on $\mathrm{int}(\mathrm{dom}(\varphi)) = \mathbb{R}^{n}$. Since $\mathcal{V}$ is strictly convex, by invoking  \cite{Schneider_2013_Convex_Body}, Corollary 1.7.3 we have $\|\cdot\|_{\mathcal{V}^{\circ}}$ is differentiable on $\mathbb{R}^{n}\setminus \{0\}$. Additionally, given $g$ is differentiable on $\mathbb{R}$, $g\circ \|\cdot\|_{\mathcal{V}^{\circ}}$ is differentiable on $\mathbb{R}^{n}\setminus \{0\}$. It remains to show $\varphi$ is differentiable at $0$, which is equivalent to verifying that the normal cone of $\mathrm{epi}(\varphi)$ at $(0,\varphi(0))$  is a single ray. Specifically, we show $\{\gamma(0,-1):\gamma\geq 0\}$ is the unique ray in $\mathcal{N}_{\mathrm{epi}(\varphi)}((0,g\circ\|0\|_{\mathcal{V}^{\circ}})) = \mathcal{N}_{\mathrm{epi}(\varphi)}((0,0)) $ (\ie, $\nabla\varphi(0) = 0 $). We begin by showing $(0,-1) \in \mathcal{N}_{\mathrm{epi}(\varphi)}((0,0))$. For any $(x,t) \in \mathrm{epi}(\varphi)$, $\langle (0,-1), (x,t) - (0,0) \rangle = -t \leq 0 $ (Given $\mathrm{epi}(\varphi) \subset \mathbb{R}^{n} \times \mathbb{R}_{+}$ ). 

Next we show that for any $(v,-1)$ such that $v\neq 0$, it holds that $(v,-1) \notin \mathcal{N}_{\mathrm{epi}(\varphi)}((0,0)) $. Let $x = \alpha v$ for some sufficiently small $\alpha > 0$, we have $(x,\varphi(x)) = (\alpha  v, g(\|\alpha v\|_{\mathcal{V}^{\circ}})) = (\alpha v,  g(\alpha\|v\|_{\mathcal{V}^{\circ}})) \in \mathrm{epi}(\varphi)$. We are interested in if the following dot product is positive $\langle (v,-1), (x,\varphi(x)) - (0,0) \rangle = \langle (v,-1), (x,\varphi(x)) \rangle=\alpha \|v\|_{2}^{2} -g(\alpha\|v\|_{\mathcal{V}^{\circ}})  $. Since $\alpha>0$, we can equivalently show $\|v\|_{2}^{2} - \frac{g(\alpha\|v\|_{\mathcal{V}^{\circ}}) }{\alpha}$ is positive for some sufficiently small $\alpha$. Indeed,
$
\lim_{\alpha \rightarrow 0^{+}} \|v\|_{2}^{2} - \frac{g(\alpha\|v\|_{\mathcal{V}^{\circ}}) }{\alpha} 
= \|v\|_{2}^{2} - \lim_{\alpha \rightarrow 0^{+}} \frac{g(\alpha\|v\|_{\mathcal{V}^{\circ}}) }{\alpha}
=\|v\|_{2}^{2} - \lim_{\gamma \rightarrow 0^{+}} \frac{g(\gamma) }{\gamma} \|v\|_{\mathcal{V}^{\circ}} = \|v\|^{2}_{2} > 0,
$
where the third equality is by the differentiability of $g$ at $0$, specifically $\nabla g(0)=0$. The inequality is due to $v\neq 0$. In conclusion, $\varphi$ is differentiable on $\mathbb{R}^{n}$, with $\nabla \varphi(0) = 0$.

\emph{(B) Legendre: Essentially smooth.} Having proved the differentiability of $\varphi$ on $\mathbb{R}^{n}$, to conclude $\varphi$ is essentially smooth, it suffice to show $\|\nabla \varphi(x_{n})\|_{2} \rightarrow + \infty $, for any sequence $(x_{n}) \subset \mathrm{int}(\mathrm{dom}(\varphi))$ such that $x_{n} \rightarrow x \in \mathrm{bd}(\mathrm{dom}(\varphi))$. Given $\mathrm{dom}(\varphi)= \mathbb{R}^{n}$, $\mathrm{bd}(\mathrm{dom}(\varphi)) = \varnothing$, the condition is vacuously satisfied.

\emph{(C) Zero at Origin.} We have proved $\varphi(0)=0$ and $\nabla \varphi (0) = 0$ in the course of proving the differentiability of $\varphi$. \hfill \Halmos
\end{proof}

Then we prove that under a Legendre distance generating function, the induced Bregman proximal operator has a Bregman projection interpretation.
\begin{claim}[Bregman Proximal Operator is Bregman Projection]
\label{claim: prox_is_proj}
Fix a function $\psi: \mathbb{R}^{n} \rightarrow \mathbb{R} \cup \{+\infty\}$ that is Legendre, and a closed convex set $\mathcal{S} \subset \mathrm{int}(\mathrm{dom}(\psi))$. Then
$$\mathrm{Prox}^{\psi}_{c,\mathcal{S}}(y,\eta) = \Pi^{\psi}_{\mathcal{S}}\left(\nabla \psi ^{*} (\nabla \psi (y) - \eta c ) \right).$$
\end{claim}
\begin{proof}{Proof.}
Denote $\overline{x}=  \mathrm{Prox}^{\psi}_{c,\mathcal{S}}(y,\eta)=\argmin_{x\in \mathcal{S}} \langle c,x \rangle + \frac{1}{\eta} D_{\psi}(x, y)$. The minimizer is unique due to strong convexity of Legendre distance-generating functions (Assumption \ref{assumption: blanket varphi}). By the first order optimality condition, $\langle \eta c + \nabla \psi(\overline{x}) - \nabla \psi(y), x' - \overline{x} \rangle \geq 0, \ \forall x' \in \mathcal{S}$, which is exactly the first order optimality condition (Lemma \ref{lemma: Bregman Projection VI}) for $\overline{x} = \Pi^{\psi}_{\mathcal{S}}\left((\nabla \psi) ^{-1} (\nabla \psi (y) - \eta c ) \right)=\Pi^{\psi}_{\mathcal{S}}\left(\nabla \psi ^{*} (\nabla \psi (y) - \eta c ) \right).$
\hfill \Halmos
\end{proof}

\begin{claim}\label{claim: robust path is proj}
The robust path $\mathcal{P}'(\mathcal{V})$ of $(\mathrm{RC})$ is equivalent to the Bregman Projection induced by $\varphi$ of curve $\left\{   \nabla \varphi^{*}\left(\nabla\varphi(0) - \omega^{-1} a_{0} \right) : \omega \in [0,\infty) \right\}$ onto the feasible region $\mathcal{X}$:
    $$\mathcal{P}'(\mathcal{V})= \left\{  \Pi^{\varphi}_{\mathcal{X}} \left( \nabla \varphi^{*}\left(\nabla\varphi(0) - \omega^{-1} a_{0} \right) \right) : \omega \in [0,\infty) \right\}$$
\end{claim}
\begin{proof}{Proof.}
This result follows directly from the previous claims.
$$
\begin{aligned}
x_{\mathrm{R}}'(\omega, \mathcal{V})=& \ \argmin_{x\in \mathcal{X}} \langle a_{0}, x \rangle + \omega \varphi(x)\\
=& \ \argmin_{x\in \mathcal{X}} \langle a_{0}, x \rangle + \omega ( \varphi(x) - \varphi(0) - \langle \nabla\varphi(0), x-0 \rangle )\\
=& \ \mathrm{Prox}^{\varphi}_{a_{0}, \mathcal{X}}\left(0, \omega^{-1}\right)\\
=& \ \Pi^{\varphi}_{\mathcal{X}} \left( \nabla \varphi^{*}\left(\nabla\varphi(0) - \omega^{-1} a_{0} \right) \right).
\end{aligned}
$$ 
The second equality is due to Claim \ref{claim: gauge_set_to_Legendre_dgf} (C). The third equality is by Definition \ref{def: Bregman Divergence}. The fourth equality is due to Claim \ref{claim: gauge_set_to_Legendre_dgf} (B) and Claim \ref{claim: prox_is_proj}. \hfill \Halmos
\end{proof}

\subsubsection{Central Path and Proximal Path: Bregman Projection and Dual Space.}\label{subsubsection-Thm1-proof-CP-PP}

\begin{claim}
\begin{itemize}
    \item[Central Path:] The {\CentralPathFull} ({\CentralPathAbb}) $\{x_{\mathrm{{\CentralPathAbb}}}(\omega): \omega \in [0,\infty) \}$ induced by the distance-generating function $\varphi(\cdot)=g \circ \|\cdot\|_{\mathcal{V}^{\circ}}$ initialized at $x_{\mathrm{R}}$ can be formulated as the Bregman projection of the curve $\left\{ \nabla \varphi^{*}\left(\nabla\varphi(x_{\mathrm{R}}) - \omega^{-1} a_{0} \right)  : \omega \in [0,\infty) \right\}$ onto the feasible region $\mathcal{X}$:
    $$
    \left\{x_{\mathrm{{\CentralPathAbb}}}(\omega) : \omega \in [0,\infty) \right\}= \left\{\Pi^{\varphi}_{\mathcal{X}} \left( \nabla \varphi^{*}\left(\nabla\varphi(x_{\mathrm{R}}) - \omega^{-1} a_{0} \right) \right) : \omega \in [0,\infty) \right\}.
    $$
    \item[Proximal Path:] The {\ProximalPathFull} ({\ProximalPathAbb}) $\{x_{k}\}$ of $(\mathrm{P})$  induced by the distance-generating function $\varphi(\cdot) = g \circ \|\cdot\|_{\mathcal{V}^{\circ}}$ initialized at $ x_{\mathrm{R}}$ and associated with sequence $\{\lambda_{k}\}$ can be formulated as successive Bregman projections onto the feasible region $\mathcal{X}$:
    $$
    x_{k+1}  = \Pi^{\varphi}_{\mathcal{X}} \left( \nabla \varphi^{*}\left(\nabla\varphi(x_{k}) - \lambda_{k}^{-1} a_{0} \right) \right) , \quad \text{for } k  = 0,1,\cdots.
    $$
\end{itemize}
\end{claim}
\begin{proof}{Proof.}
The proof is a straightforward application of Definition \ref{def: proximal operator} and Claim \ref{claim: prox_is_proj}.
$$
\begin{aligned}
x_{\mathrm{{\CentralPathAbb}}}(\omega) =& \ \argmin_{x\in \mathcal{X}} \langle a_0,x  \rangle + \omega D_{\psi}(x,x_{0})\\
=& \ \mathrm{Prox}^{\psi}_{a_{0},\mathcal{X}}(x_{0},\omega^{-1})\\
=& \ \Pi^{\varphi}_{\mathcal{X}} \left( \nabla \varphi^{*}\left(\nabla\varphi(x_{\mathrm{R}}) - \omega^{-1} a_{0} \right) \right). 
\end{aligned}
$$
$$
\begin{aligned}
x_{k+1} =& \ \argmin_{x\in \mathcal{X}} \langle a_0,x  \rangle + \lambda_{k}D_{\psi}(x,x_{k}) \\ 
=& \ \mathrm{Prox}^{\psi}_{a_{0},\mathcal{X}}(x_{k},\lambda_{k}^{-1}) \\
=& \ \Pi^{\varphi}_{\mathcal{X}} \left( \nabla \varphi^{*}\left(\nabla\varphi(x_{k}) - \lambda_{k}^{-1} a_{0} \right) \right).
\end{aligned}
$$
\hfill \Halmos
\end{proof}

\section{Recovering Robust Path by Optimization Paths}
\label{section: approxiamte_robust_path}

In this section, leveraging the geometric view of Section \ref{section: RP_geometry}, we show that the robust path $\mathcal{P}'(\mathcal{V})$ of $(\mathrm{RC})$ can be approximated by some appropriately designed {\CentralPathFull} $\{x_{\mathrm{{\CentralPathAbb}}}(\omega)\}$ and {\ProximalPathFull} $\{x_{k}\}$ of $(\mathrm{P})$. The two optimization paths' distance-generating functions $\varphi$ are induced by the uncertainty set shape $\mathcal{V}$, \ie, $\varphi(\cdot) = g \circ \|\cdot\|_{\mathcal{V}^{\circ}}$. Furthermore, the two optimization paths are initialized at the most robust solution $x_{\mathrm{R}}$, while solving for the most efficient solution $x_{\mathrm{E}}$. 


The structure of the main results in this section is summarized in Figure \ref{fig: summary_sec_2_results}. In Section \ref{subsec: BPP is RP}, we establish that {\CentralPathFull}s are good approximations of, and sometimes exact, robust paths (Theorem \ref{theorem: uniform bound CP and RP} and Corollary \ref{corollary: BPP is RP}). This builds on the previous insight that the robust path and the central path are both Bregman projections of similar curves onto the feasible region, where they differ only in the initial points of the curves (Theorem \ref{theorem: RP_CP_PP are Bregman Projections}).

Theorem \ref{theorem: RP_CP_PP are Bregman Projections} also points out that the central path and the proximal path share the same initial point, but differ only in that the central path is the Bregman projection of an entire curve onto the feasible region, while the proximal path is generated via \emph{successive} Bregman projections of short curve segments onto the same feasible region. Building on this insight, in Section \ref{subsec: PPM is BPP}, we first show two special cases where {\CentralPathFull}s are exact {\ProximalPathFull}s (Propositions \ref{prop: monotone_BPP_is_PPM} and \ref{prop: X_matches_V}). We then build on this analysis to establish a general bound between {\CentralPathFull}s and {\ProximalPathFull}s (Theorem \ref{theorem: PPM as approx BPP}). 

In Section \ref{sec: PPMP_Exact_RP}, we close the loop and give sufficient conditions under which {\ProximalPathFull}s are exact robust paths (Theorem \ref{theorem: PPMP_Exact_RP}). This theorem also quantifies the exact relationship between the step-size progression on a proximal path and the uncertainty set radii on the robust path.

Finally, inspired by this theory, we state an algorithm for (approximately) recovering a robust path via a single algorithmic pass of the proximal method. Numerical validations in the subsequent Section \ref{section: application} show that this algorithm verifies the predictions of our theory when the problem's technical setup matches our theory. Surprisingly, we also observe that this algorithm retains a strong performance (in terms of solution quality and computational time) even when the technical conditions are severely violated.



\subsection{{Robust Path and \CentralPathFullCap}}
\label{subsec: BPP is RP}
We define the $\kappa$-expansiveness of a Bregman projection operator as follows.
\begin{definition}
\label{def: Kappa_expansive}
    Give a distance-generating function, $\varphi$, the induced Bregman projection $\Pi^{\varphi}_{\mathcal{S}}: \mathrm{int}(\mathrm{dom}(\varphi)) \rightarrow \mathcal{S}$ is $\kappa$-expansive if for any closed and convex $\mathcal{S}$, the following inequalities hold:
   $$D_{\varphi}\left(\Pi^{\varphi}_{\mathcal{S}}(x), \Pi^{\varphi}_{\mathcal{S}}(y) \right) \leq  \kappa \cdot D_{\varphi} \left(x,y \right), \ \forall x,y \in\mathrm{int}(\mathrm{dom}(\varphi)). $$
   In addition, denote $d = y-x$: $$D_{\varphi}\left(\Pi^{\varphi}_{\mathcal{S}+d}(x), \Pi^{\varphi}_{\mathcal{S}}(x) \right) \leq  \kappa \cdot D_{\varphi} \left(x,y \right), \ \forall x,y \in\mathrm{int}(\mathrm{dom}(\varphi)). $$
\end{definition}
For instance, for $\varphi = \frac{1}{2}\|x\|^{2}_{2}$, $\Pi^{\varphi}_{\mathcal{S}}$  is the usual Euclidean projection and is 1-expansive. More information can be found in Appendix \ref{sec:kappa-expansiveness-example}.

We now show that the robust path $\mathcal{P}'(\mathcal{V})$ of $(\mathrm{RC})$  can be approximated by the {\CentralPathFull} $\{x_{\mathrm{{\CentralPathAbb}}}(\omega)\}$ for $(\mathrm{P})$ induced by $\varphi(\cdot) = g\circ \|\cdot\|_{\mathcal{V}^{\circ}}$ and initialized at $x_{\mathrm{R}}$. Its proof relies on the geometric view on both $\mathcal{P}'(\mathcal{V})$ and $\{x_{\mathrm{{\CentralPathAbb}}}(\omega)\}$ made possible by Theorem \ref{theorem: RP_CP_PP are Bregman Projections}. 
\begin{theorem}
\label{theorem: uniform bound CP and RP}
Assume $\mathcal{V}$ satisfies Assumption \ref{assumption: blanket V}, $\varphi$ satisfies Assumption \ref{assumption: blanket varphi}, and the induced Bregman projection $\Pi^{\varphi}_{\mathcal{S}}$ is $\kappa$-expansive. The Bregman divergence between the central path and the robust path is uniformly bounded by an upper bound. The upper bound depends on $\kappa$ and the Bregman divergence between two points.
$$
\begin{aligned}
D_{\varphi}\left( x_{\mathrm{{\CentralPathAbb}}}(\omega), x_{\mathrm{R}}'(\omega) \right) &\leq \kappa^{2} \cdot D_{\varphi}\left( \Pi^{\varphi}_{\mathcal{X}} (0), \Pi^{\varphi}_{\mathrm{Aff}(\mathcal{X})} (0) \right) \\
\text{and \quad} D_{\varphi}\left(x_{\mathrm{R}}'(\omega), x_{\mathrm{{\CentralPathAbb}}}(\omega)  \right) &\leq \kappa^{2} \cdot D_{\varphi}\left(\Pi^{\varphi}_{\mathrm{Aff}(\mathcal{X})} (0) , \Pi^{\varphi}_{\mathcal{X}} (0)\right), \quad  \forall   \omega \in [0,\infty).  
\end{aligned}
$$
Moreover, the bound is sharp. 
\end{theorem}

\begin{proof}{Proof.}
On a high-level, the proof proceeds by first mapping the objects of interest via bijection $\nabla \varphi: \mathrm{int}(\mathrm{dom}(\varphi)) \rightarrow \mathrm{int}(\mathrm{dom}(\varphi^{*}))$ from the primal space to the dual space where the paths enjoy simple structure, before mapping the objects via $(\nabla \varphi)^{-1} = \nabla \varphi^{*}: \mathrm{int}(\mathrm{dom}(\varphi^{*})) \rightarrow \mathrm{int}(\mathrm{dom}(\varphi))$ back to the primal space to establish the result.

We begin by taking the geometry view of the two paths of Theorem \ref{theorem: RP_CP_PP are Bregman Projections}. Together with Lemma \ref{lemma: Breg_Proj_calc_1}, we have
\begin{equation}
\label{eq: uniform bound 1}
        x_{\mathrm{R}}'(\omega) 
        =  \Pi^{\varphi}_{\mathcal{X}} \circ \Pi^{\varphi}_{\mathrm{Aff}(\mathcal{X})}\left(\nabla \varphi ^{*} \left(\nabla \varphi (0) - \omega^{-1}a_{0} \right) \right),
\end{equation}
\begin{equation}
\label{eq: uniform bound 2}
        x_{\mathrm{{\CentralPathAbb}}}(\omega)
        = \Pi^{\varphi}_{\mathcal{X}}\circ \Pi^{\varphi}_{\mathrm{Aff}(\mathcal{X})}\left(\nabla \varphi ^{*} \left(\nabla \varphi (x_{\mathrm{R}}) - \omega^{-1}a_{0} \right) \right).
\end{equation}
We recall by Definition \ref{definition: unique robust path}, $x_{\mathrm{R}} = \argmin_{x\in \mathcal{X}} \varphi(x)$; together with Claim \ref{claim: gauge_set_to_Legendre_dgf}.C, we have $x_{\mathrm{R}} = \Pi^{\varphi}_{\mathcal{X}}(0)$.

Define $x_{\mathrm{A}} := \Pi^{\varphi}_{\mathrm{Aff}(\mathcal{X})} (0) \in \mathrm{Aff}(\mathcal{X})$, we have $\mathrm{Aff}(\mathcal{X}) = \mathcal{L} + x_{\mathrm{A}}$, where $\mathcal{L} = \mathrm{Aff}(\mathcal{X}) - x_{\mathrm{A}}$ is a linear subspace. Now, we apply the dual perspective of Lemma \ref{lemma:dual_cone bregman_proj} on $x_{\mathrm{A}} = \Pi^{\varphi}_{\mathrm{Aff}(\mathcal{X})} (0) $: 
\begin{subequations}
    \begin{align}
        \nabla \varphi \left( x_{\mathrm{A}}\right) =& \ \nabla \varphi \left(\Pi^{\varphi}_{\mathrm{Aff}(\mathcal{X})} (0) \right) \\
        =&  \ \nabla \varphi \left(\Pi^{\varphi}_{\mathcal{L} + \Pi^{\varphi}_{\mathrm{Aff}(\mathcal{X})} (0)} (0) \right) \\
        = & \ \Pi^{\varphi^{*}}_{\mathcal{L}^{\bot} + \nabla \varphi(0)} \left(\nabla \varphi \left(\Pi^{\varphi}_{\mathrm{Aff}(\mathcal{X})} (0) \right)\right)  \\
        = & \ \Pi^{\varphi^{*}}_{\mathcal{L}^{\bot} + \nabla \varphi(0)} \left(\nabla \varphi \left(x_{\mathrm{A}}\right)\right).    
    \end{align}
\end{subequations}
Hence $\nabla \varphi \left(x_{\mathrm{A}} \right)  \in \mathcal{L}^{\bot} + \nabla \varphi(0)$, or equivalently 
$\mathcal{L}^{\bot} + \nabla \varphi(0) = \mathcal{L}^{\bot} + \nabla \varphi \left(x_{\mathrm{A}} \right)$.

Next, we map $\Pi^{\varphi}_{\mathrm{Aff}(\mathcal{X})}\left(\nabla \varphi ^{*} \left(\nabla \varphi (0) - \omega^{-1}a_{0} \right) \right)$ and $\Pi^{\varphi}_{\mathrm{Aff}(\mathcal{X})}\left(\nabla \varphi ^{*} \left(\nabla \varphi (x_{\mathrm{R}}) - \omega^{-1}a_{0} \right) \right)$ of (\ref{eq: uniform bound 1}) and (\ref{eq: uniform bound 2}) to the dual space where the pair enjoy a simpler characterization of Bregman projection of $x_{\mathrm{A}}$ onto two parallel affine subspaces.

We begin by applying the dual characterization of Lemma \ref{lemma:dual_cone bregman_proj} on $\Pi^{\varphi}_{\mathrm{Aff}(\mathcal{X})}\left(\nabla \varphi ^{*} \left(\nabla \varphi (0) - \omega^{-1}a_{0} \right) \right)$. For any $\omega \in [0,\infty)$ we have 
\begin{subequations}
\label{eq: uniform bound 3}
    \begin{align}
         &\nabla \varphi \left( \Pi^{\varphi}_{\mathrm{Aff}(\mathcal{X})}\left(\nabla \varphi ^{*} \left(\nabla \varphi (0) - \omega^{-1}a_{0} \right) \right)\right)\\
        = &  \ \nabla \varphi \left( \Pi^{\varphi}_{\mathcal{L} + x_{\mathrm{A}}}\left(\nabla \varphi ^{*} \left(\nabla \varphi (0) - \omega^{-1}a_{0} \right) \right)\right) 
        \\= & \ \Pi^{\varphi^{*}}_{\mathcal{L}^{\bot} + \nabla \varphi (0) - \omega^{-1}a_{0}} \left( \nabla \varphi \left( x_{\mathrm{A}} \right) \right)
        \\= & \ \Pi^{\varphi^{*}}_{\mathcal{L}^{\bot} + \nabla \varphi \left(x_{\mathrm{A}}\right) - \omega^{-1}a_{0}} \left( \nabla \varphi \left(x_{\mathrm{A}}\right) \right),
    \end{align}
\end{subequations}
where the first equality is due to $x_{\mathrm{A}} \in \mathrm{Aff}(\mathcal{X})$, the second equality is a direct application of Lemma \ref{lemma:dual_cone bregman_proj}, the third equality is by $\mathcal{L}^{\bot} + \nabla \varphi(0) = \mathcal{L}^{\bot} + \nabla \varphi \left(x_{\mathrm{A}} \right)$. 

Similarly, for $\Pi^{\varphi}_{\mathrm{Aff}(\mathcal{X})}\left(\nabla \varphi ^{*} \left(\nabla \varphi (x_{\mathrm{R}}) - \omega^{-1}a_{0} \right) \right)$, we have
\begin{subequations}
\label{eq: uniform bound 4}
    \begin{align}
        &\nabla \varphi \left( \Pi^{\varphi}_{\mathrm{Aff}(\mathcal{X})}\left(\nabla \varphi ^{*} \left(\nabla \varphi (x_{\mathrm{R}}) - \omega^{-1}a_{0} \right) \right)\right) \\
        = & \ \nabla \varphi \left( \Pi^{\varphi}_{\mathcal{L} + x_{\mathrm{A}}}\left(\nabla \varphi ^{*} \left(\nabla \varphi (x_{\mathrm{R}}) - \omega^{-1}a_{0} \right) \right)\right) \\
        = & \ \Pi^{\varphi^{*}}_{\mathcal{L}^{\bot} + \nabla \varphi (x_{\mathrm{R}}) - \omega^{-1}a_{0}} \left( \nabla \varphi \left( x_{\mathrm{A}} \right) \right).
    \end{align}
\end{subequations}
Now we are ready to prove the Theorem:
$$
    \begin{aligned}
        & D_{\varphi}\left( x_{\mathrm{{\CentralPathAbb}}}(\omega), x_{\mathrm{R}}'(\omega) \right) \\
        = & \  D_{\varphi}\left( \Pi^{\varphi}_{\mathcal{X}} \circ \Pi^{\varphi}_{\mathrm{Aff}(\mathcal{X})}\left(\nabla \varphi ^{*} \left(\nabla \varphi (0) - \omega^{-1}a_{0} \right) \right),  \Pi^{\varphi}_{\mathcal{X}} \circ \Pi^{\varphi}_{\mathrm{Aff}(\mathcal{X})}\left(\nabla \varphi ^{*} \left(\nabla \varphi (x_{\mathrm{R}}) - \omega^{-1}a_{0} \right) \right) \right) \\
        \leq & \ \kappa  D_{\varphi}\left(  \Pi^{\varphi}_{\mathrm{Aff}(\mathcal{X})}\left(\nabla \varphi ^{*} \left(\nabla \varphi (0) - \omega^{-1}a_{0} \right) \right),  \Pi^{\varphi}_{\mathrm{Aff}(\mathcal{X})}\left(\nabla \varphi ^{*} \left(\nabla \varphi (x_{\mathrm{R}}) - \omega^{-1}a_{0} \right) \right) \right) \\
        =& \ \kappa  D_{\varphi^{*}}\left( \nabla \varphi \left( \Pi^{\varphi}_{\mathrm{Aff}(\mathcal{X})}\left(\nabla \varphi ^{*} \left(\nabla \varphi (x_{\mathrm{R}}) - \omega^{-1}a_{0} \right) \right)\right),  \nabla \varphi \left( \Pi^{\varphi}_{\mathrm{Aff}(\mathcal{X})}\left(\nabla \varphi ^{*} \left(\nabla \varphi (0) - \omega^{-1}a_{0} \right) \right)\right)\right) \\
        =& \ \kappa  D_{\varphi^{*}}\left(\Pi^{\varphi^{*}}_{\mathcal{L}^{\bot} + \nabla \varphi (x_{\mathrm{R}}) - \omega^{-1}a_{0}} \left( \nabla \varphi \left( x_{\mathrm{A}} \right) \right),  \Pi^{\varphi^{*}}_{\mathcal{L}^{\bot} + \nabla \varphi \left(x_{\mathrm{A}}\right) - \omega^{-1}a_{0}} \left( \nabla \varphi \left(x_{\mathrm{A}}\right) \right) \right) \\
        \leq & \ \kappa^{2} D_{\varphi^{*}}(\nabla \varphi \left(x_{\mathrm{A}}\right), \nabla \varphi \left(x_{\mathrm{R}}\right)) \\
        =& \ \kappa^{2}D_{\varphi}\left(x_{\mathrm{R}},x_{\mathrm{A}} \right).
    \end{aligned}
$$
The first equality follows from (\ref{eq: uniform bound 1}) and (\ref{eq: uniform bound 2}), the first inequality is by Definition \ref{def: Kappa_expansive}, the second equality is by \cite{bauschke1997legendre} Theorem 3.7(v) (mapping from the primal space to the dual space), the third equality is by (\ref{eq: uniform bound 3}) and (\ref{eq: uniform bound 4}), the second inequality is by Definition \ref{def: Kappa_expansive}, the fourth equality is again by \cite{bauschke1997legendre} Theorem 3.7(v) (mapping from the dual space back to the primal space). The other inequality $D_{\varphi}\left(x_{\mathrm{R}}'(\omega), x_{\mathrm{{\CentralPathAbb}}}(\omega)  \right) \leq \kappa^{2} \cdot D_{\varphi}\left(x_{\mathrm{A}},x_{\mathrm{R}}\right)$ follows from the same argument.

To finish the proof, we establish the sharpness of the above bound via the following example:
\begin{example}[Upper Bound Sharpness] Consider following problem instance
    $a_{0}=(-1,1)$, $\mathcal{X}=\{(x_{1}, x_{2}) \in \mathbb{R}^{2}: x_{1} + 2x_{2} = 2, 
 \ x_{1}\geq 0.5, \ x_{2} \geq 0 \}$, with $\ell_2$ norm setup, $\varphi(x) = \frac{1}{2}\|x\|^{2}_{2}$. It can be easily  verified that for $\omega=0.5$, $\frac{1}{2}\|x_{\mathrm{CP}}(\omega)-x_{\mathrm{R}}'(\omega)\|^{2}_{2} = \frac{1}{2}\|  \Pi_{\mathcal{X}} (0)-\Pi_{\mathrm{Aff}(\mathcal{X})} (0)\|^{2}_{2}$ or equivalently, $D_{\varphi}\left( x_{\mathrm{CP}}(\omega), x_{\mathrm{R}}'(\omega) \right) = \kappa^{2} D_{\varphi}\left( \Pi^{\varphi}_{\mathcal{X}} (0), \Pi^{\varphi}_{\mathrm{Aff}(\mathcal{X})} (0) \right)$, where $\kappa = 1$. \hfill \Halmos
\end{example}
\end{proof}

\begin{corollary}
    \label{corollary: BPP is RP}
Assume $\mathcal{V}$ satisfies Assumption \ref{assumption: blanket V}, and $\varphi$ satisfies Assumption \ref{assumption: blanket varphi}, then the following statement is true.
$$\Pi^{\varphi}_{\mathcal{X}} (0) = \Pi^{\varphi}_{\mathrm{Aff}(\mathcal{X})} (0) \quad \Longrightarrow \quad x_{\mathrm{R}}'(\omega)= x_{\mathrm{{\CentralPathAbb}}}(\omega), \quad  \forall  \omega \in [0,\infty).$$
\end{corollary}
\begin{proof}{Proof.}
    This result follows directly from Theorem \ref{theorem: uniform bound CP and RP}: the right hand side of Theorem \ref{theorem: uniform bound CP and RP}'s inequality collapses to zero when the two initial points match each other, $\Pi^{\varphi}_{\mathcal{X}} (0) = \Pi^{\varphi}_{\mathrm{Aff}(\mathcal{X})} (0)$. \Halmos
\end{proof}

\begin{example}[Zero Gap Between Robust Path and Central Path]
    $\Pi^{\varphi}_{\mathcal{X}} (0) = \Pi^{\varphi}_{\mathrm{Aff}(\mathcal{X})} (0)$ for the following cases.
    \begin{itemize}
        \item $\mathcal{X}$ contains 0.
        \item $\Pi^{\varphi}$ is the Euclidean projection, and $\mathcal{X}$ is the intersection of a positively oriented affine space with the positive orthant: $\mathcal{X} = \{ x : Ax = b, x \geq 0 \}$ where $A^\intercal (AA^\intercal)^{-1}b \geq 0$. For instance, a simplex $\mathcal{X} = \{x\in \mathbb{R}^{n}_{+}: \langle 1,x\rangle = 1\}$ satisfies this criterion.
    \end{itemize}
\end{example}

\subsection{{\CentralPathFullCap} and {\ProximalPathFullCap}}
\label{subsec: PPM is BPP}
In this section, we show the two optimization paths of $(\mathrm{P})$, \ie, the {\ProximalPathFull} $\{x_k\}$ and the {\CentralPathFull} $\{x_{\mathrm{{\CentralPathAbb}}}(\omega)\}$ are in general approximations of each other, and under two special cases equivalent. 
\subsubsection{Special Case One: Polyhedral Monotonicity of Feasible Regions.}

For convex polyhedron feasible regions $\mathcal{X}$, we show that if the following is true: once $\{x_{k}\}$ enters a face of $\mathcal{X}$ it remains in that face, then $\{x_{k}\} = \{x_{\mathrm{{\CentralPathAbb}}}(\omega_{k})\}$. More precisely, we adopt the following definition from \cite{monotone_QRLP}.
\begin{definition}
\label{def:monotone}
    The {\ProximalPathFull} $\{x_{k}\}$ is \textit{monotone} on convex polyhedron $\mathcal{X}$, if for any face $\mathcal{F}$ of $\mathcal{X}$: $\quad  x_{k} \in \mathcal{F} \quad  \Rightarrow \quad  x_{k+n} \in \mathcal{F}, \ \forall n\in [1,K].$
\end{definition}
Intuitively, monotonicity can be interpreted as: once a polyhedron constraint becomes active for $x_{k}$, it stays active for all the subsequent sequences. Another equivalent interpretation is that $\{x_{k}\}$ generates a path in the partially ordered set (poset) of the faces of $\mathcal{X}$ that is nonincreasing in the set order.

\begin{proposition}
\label{prop: monotone_BPP_is_PPM}
Let the feasible region $\mathcal{X}$ be a convex polyhedron. Let $\{x_k\}$ be a {\ProximalPathFull} initialized by $x_0$ associated with the step-size sequence $\{\lambda_{k}> 0 : \, \sum_{k=0}^{\infty} \lambda_{k}^{-1} = \infty \}$. Let $\{x_{\mathrm{{\CentralPathAbb}}}(\omega_{k})\}$ be the {\CentralPathFull} initialized also at $x_0$. If any of the following conditions is satisfied:
\begin{itemize}
    \item[(C1. Monotone):] $\{x_{k}\}$ is monotone on $\mathcal{X}$,
    \item[(C2. Affine Subspace):] $\mathcal{X}= \{x \in \mathbb{R}^{n}: Ax = b \}$,
    \item[(C3. Unconstrained):] $\mathcal{X} = \mathbb{R}^{n}$,
\end{itemize}
then
$x_{k+1} = x_{\mathrm{{\CentralPathAbb}}}\left(\omega_{k+1}\right), \ \forall k\in [0,K].$ Furthermore, $\omega_{k}$ can be recovered in closed form as a function of the proximal path step-size sequence $\{\lambda_{k}\}$:
 $\omega_{k} = \left( \lambda_{0}^{-1} + \dots + \lambda_{k-1}^{-1} \right)^{-1}.$
\end{proposition}

\begin{figure}[t]
    \centering
    \begin{subfigure}[t]{0.43\linewidth}
        \centering
        \includegraphics[width=\linewidth]{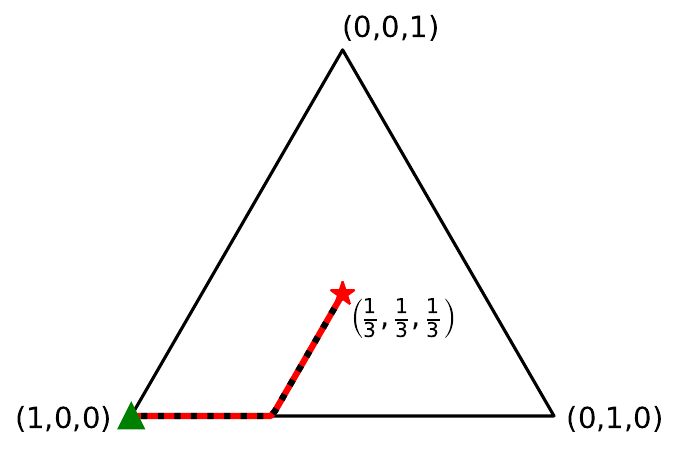}
        \caption{\footnotesize $\{x_{k}\}$ is monotone on $\mathcal{X}=\Delta^{3} \\  \xRightarrow{\text{Prop. \ref{prop: monotone_BPP_is_PPM}}} \{x_{k}\} = \{x_{\mathrm{{\CentralPathAbb}}}\left(\omega_{k}\right)\}.$}
        \label{fig:PPM_is_BPP_A}
    \end{subfigure}%
    \begin{subfigure}[t]{0.46\linewidth}
        \centering
        \includegraphics[width=\linewidth]{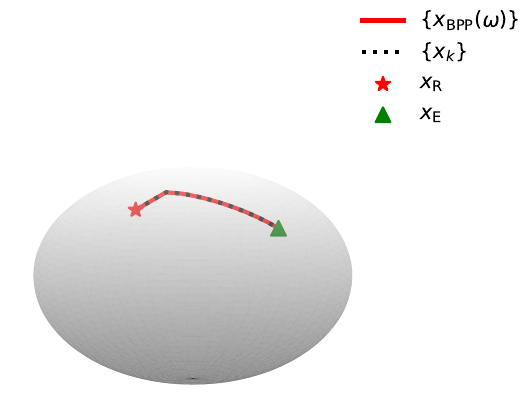}
        \caption{\footnotesize $\mathcal{U} = \{x\in \mathbb{R}^{3}: \langle (x-c),A(x-c)\rangle \leq 1 \}$, the feasible region is $\mathcal{U}$, the uncertainty set is $\mathcal{U}^{\circ}$ $\xRightarrow{\text{Prop. \ref{prop: X_matches_V}}} \{x_k\} \subset \{x_{\mathrm{{\CentralPathAbb}}}\left(\omega\right): \omega \in [0,\infty)\}$.}
        \label{fig:PPM_is_BPP_B}
    \end{subfigure}%
    \caption{Examples of equivalence results between {\ProximalPathFull} $\{x_{k}\}$ and {\CentralPathFull} $\{x_{\mathrm{{\CentralPathAbb}}}\left(\omega\right)\}$.}
    \label{fig:PPM_is_BPP}
\end{figure}

\begin{proof}{Proof.}
We defer the proof to Appendix \ref{sec: proof prop: monotone_BPP_is_PPM}.

\end{proof}
\textit{Remark}. The monotonicity condition (C1) can be verified after generating the {\ProximalPathFull} $\{x_k\}$, \ie, verify if for all inequality constraints, once a constraint is active for $x_k$, it remains active for $x_{k+n}$ for all $n\in[1,K]$. A canonical example that satisfies the monotonicity condition (C1) is a simplex feasible region $\mathcal{X} = \Delta^{n} = \{x\in \mathbb{R}^{n}_{+}: \langle 1,x\rangle = 1\}$ under an ellipsoidal uncertainty set  $\mathcal{V}$ as depicted in Figure \ref{fig:PPM_is_BPP_A}.

\subsubsection{Special Case Two: Feasible Region and Uncertainty Set are Polar Pairs.}
Consider the following motivating example: if the feasible region $\mathcal{X}$ and the uncertainty set $\mathcal{V}$ are both $2$-norm balls, then $\{x_{k}\}$ is not monotone on $\mathcal{X}$ while $\{x_k\} \subset \{x_{\mathrm{{\CentralPathAbb}}}\left(\omega\right): \omega \in [0,\infty)\}$. Most interestingly, this observation generalizes: assume the feasible region $\mathcal{X}$ can be defined as the sub-level set of a gauge function $\|\cdot\|_{\mathcal{U}} $, if the uncertainty set is chosen as $\mathcal{V}=\mathcal{U}^{\circ}$, \ie the feasible region and the uncertainty set are polar pairs up to a rescaling, then we have $\{x_k\} \subset \{x_{\mathrm{{\CentralPathAbb}}}\left(\omega\right): \omega \in [0,\infty)\}$. Such an example is shown in Figure \ref{fig:PPM_is_BPP_B}. Formally, we have the following result.

\begin{proposition}
    \label{prop: X_matches_V}
    Assume $\mathcal{X}=\{ x\in \mathbb{R}^{n}: \ \|x\|_{\mathcal{U}} \leq l\}$ where $\mathcal{U}$ satisfies Assumption \ref{assumption: blanket V}, if robust path uncertainty set is designed to be $\mathcal{U}^{\circ}$ and the corresponding optimization path distance-generating function $\varphi(\cdot) = g \circ \|\cdot\|_{\mathcal{U}}$ satisfies Assumption \ref{assumption: blanket varphi}, then: $\{x_k\} \subset \{x_{\mathrm{{\CentralPathAbb}}}\left(\omega\right): \omega \in [0,\infty)\}.$
\end{proposition}
\begin{proof}{Proof.}
We defer the proof to Appendix \ref{sec: proof prop: X_matches_V}.
\end{proof}

\subsubsection{General Case.}
In general, the {\ProximalPathFull} and the {\CentralPathFull} do not coincide. To this end, we establish a theoretical characterization of the distance between the two algorithmic paths with the next result. Reaching the result relies solely on our established results of Corollary \ref{corollary: BPP is RP} and Proposition \ref{prop: monotone_BPP_is_PPM}. Intuitively, we can partition the pair of paths by the sequence of minimal faces they traverse and analyze the distance between the pair locally on each minimal face.

\begin{theorem}
\label{theorem: PPM as approx BPP}
    Under Assumptions 1 and 2, denote $\{x_k\} = \{x_k: k\in [0,K] \}$ as the {\ProximalPathFull} induced by $\varphi$ whose Bregman projection is $\kappa$-expansive, initialized at $x_{0}$ and associated with a step-size sequence $\{\lambda_{k}> 0 : \, \sum_{k=0}^{\infty} \lambda_{k}^{-1} = \infty \}$. Denote $\upsilon = \omega^{-1}$ and let $\{x_{\mathrm{{\CentralPathAbb}}}(\upsilon^{-1}; x_{0})\}=\{x_{\mathrm{{\CentralPathAbb}}}(\upsilon^{-1}; x_{0}): \upsilon \in [0,\infty)\}$ be the {\CentralPathFull} induced by $\varphi$ and initialized at ${x_{0}}$. We assume $\{x_k\}, \{x_{\mathrm{{\CentralPathAbb}}}(\upsilon^{-1}; x_{0})\} \subset \bigcup_{i\in [I]} \mathrm{ri} \left(\mathcal{F}_{i} \right)$, where $\mathcal{F}_{i}, \ \forall i\in [I]$ are faces of $\mathcal{X}$. For every $i\in [I]$, define $ \mathcal{K}^{(i)} = \left[\underline{k}^{(i)}, \overline{k}^{(i)}\right]=  \left\{k\in [0,K]: x_{k} \in \mathrm{ri} \left(\mathcal{F}_{i} \right) \right\}$ and $ \Upsilon^{(i)} = \left[\underline{\upsilon}^{(i)}, \overline{\upsilon}^{(i)}\right]=  \left\{\upsilon\in [0,\infty): x_{\mathrm{{\CentralPathAbb}}}(\upsilon^{-1}; x_{0}) \in \mathrm{ri} \left(\mathcal{F}_{i} \right) \right\}$, consequently $\left\{x_{k}: k\in \mathcal{K}^{(i)} \right\}_{i\in [I]}$ and $\left \{x_{\mathrm{{\CentralPathAbb}}}(\upsilon^{-1}; x_{0}):  \upsilon \in \Upsilon^{(i)} \right\}_{i\in [I]}$ form partitions of $\{x_k\}$ and $\{x_{\mathrm{{\CentralPathAbb}}}(\upsilon^{-1}; x_{0})\}$ respectively. Then, for each $i\in [I]$ we have
    $$D_{\varphi} \left(x_{k}, x_{\mathrm{{\CentralPathAbb}}}(\upsilon_k^{-1}; x_{0}) \right) \leq \kappa \cdot D_{\varphi} \left( x_{\mathrm{{\CentralPathAbb}}}\left((\underline{\upsilon}^{(i)})^{-1}; x_{0} \right) , x_{\underline{k}^{(i)}}\right), \ \quad \forall k\in \left[\underline{k}^{(i)}+1, \overline{k}^{(i)}\right], $$
    where $\upsilon_{k} = \underline{\upsilon}^{(i)} +  \sum_{j =\underline{k}^{(i)}}^{k-1} \lambda_{j}^{-1} $.
\end{theorem}

\begin{proof}{Proof.}
We defer the proof to Appendix \ref{sec: proof theorem: PPM as approx BPP}.
\end{proof}

\subsection{Sufficient Condition: {\ProximalPathFullCap}s are Exact Robust Paths}
\label{sec: PPMP_Exact_RP}
In this section, we close the loop and give a sufficient condition for the {\ProximalPathFull} to be an exact robust path. 
\begin{theorem}
\label{theorem: PPMP_Exact_RP}    
Let $\{x_{k}\}$ be the {\ProximalPathFull} for $\mathrm{(P)}$ induced by $\varphi(\cdot)=g\circ \|\cdot\|_{\mathcal{V}^{\circ}}$, initialized at $x_{\mathrm{R}}$ and associated with $\{\lambda_{k}> 0 : \, \sum_{k=0}^{\infty} \lambda_{k}^{-1} = \infty \}$. Assume $\Pi^{\varphi}_{\mathcal{X}}(0) = \Pi^{\varphi}_{\mathrm{Aff}(\mathcal{X})}(0)$ and $\{x_{k}\}$ is monotone on $\mathcal{X}$, then for every $k$, $x_{k}$ is a solution to $(\mathrm{RC})$:
$$x_{k} \in \argmin_{x\in X} \max_{\xi \in \Xi(r_{k}, \mathcal{V})} \langle a_{0} + \xi, x\rangle,$$
where the corresponding uncertainty set radius $r_{k}$ admits the following closed-form expression:
$r_{k} = \omega_{k}  \nabla g \left(\|x_{k}\|_{\mathcal{V}^{\circ}}\right),$
with $\omega_{k} = \left( \lambda_{0}^{-1} + \dots + \lambda_{k-1}^{-1} \right)^{-1}.$
\end{theorem}
\begin{proof}{Proof.}
The result is a direct consequence of Corollary \ref{corollary: BPP is RP} and Proposition \ref{prop: monotone_BPP_is_PPM}:
$$
\begin{aligned}
    x_{k} \overset{\text{(Prop. \ref{prop: monotone_BPP_is_PPM}})}{=} &  x_{\mathrm{{\CentralPathAbb}}}(\omega_{k}) \\ \overset{\text{(Cor. \ref{corollary: BPP is RP}})}{=} & x_{\mathrm{R}}'(\omega_{k}) \\ \overset{\text{(Lem. \ref{lemma: unique robust path}})}{\in} & \argmin_{x\in \mathcal{X}} \max_{\xi \in \Xi(r(\omega_{k}), \mathcal{V})} \langle a_{0} + \xi, x\rangle,  
\end{aligned} 
$$ 
where $r_{k} = \omega_{k}  \nabla g \left(\|x_{k}\|_{\mathcal{V}^{\circ}}\right)$ by Lemma \ref{lemma: unique robust path}, and $\omega_{k} = \left( \lambda_{0}^{-1} + \dots + \lambda_{k-1}^{-1} \right)^{-1}$ by Proposition \ref{prop: monotone_BPP_is_PPM}. \hfill \Halmos
\end{proof}
 \textit{Remark.} Theorem \ref{theorem: PPMP_Exact_RP} gives the following algorithmic insight for generating robust paths of (RC) via {\ProximalPathFull}s of (P): (1) the design of the robust path uncertainty set shape, $\mathcal{V}$ is equivalent to the choice of the {\ProximalPathFull} distance-generating function, $\varphi$; (2) adjusting the cadence of the robust solutions' radii, $r$ corresponds to adjusting the step-size of the {\ProximalPathFull}, Theorem \ref{theorem: PPMP_Exact_RP} provides a closed-form expression of the robust solutions' $r$ as a function of the {\ProximalPathFull} solutions and the step-size sequence. 

\subsection{Algorithm: Recovering Robust Path Approximately via  {\ProximalPathFullCap}}
Theorem \ref{theorem: PPMP_Exact_RP} directly points to Algorithm \ref{alg: RP_via_PP} for recovering an (approximate) robust path of $(\mathrm{RC})$ via a single proximal path of $(\mathrm{P})$: select an appropriate step-size sequence $\{\lambda_{k}\}$ and construct an proximal path distance generating function from the robust path uncertainty set, i.e., $\varphi(\cdot) = g \circ \|\cdot\|_{\mathcal{V}^{\circ}}: \mathcal{X} \rightarrow \mathbb{R}$ satisfying Assumption  \ref{assumption: blanket varphi}; solve for the robust solution $x_{\mathrm{R}}$; initialized at $x_{\mathrm{R}}$, generate a proximal path of ($\mathrm{P}$) associated with $\{\lambda_{k}\}$ and $\varphi$; the resulting proximal path is an (approximate) robust path of ($\mathrm{RC}$), where the corresponding uncertainty set radius is a closed-form function of the proximal path solutions and the step-size sequence.

\begin{algorithm}[ht!]
\caption{Recovering (Approximate) Robust Path of $(\mathrm{RC})$ via {\ProximalPathFullCap} of $(\mathrm{P})$} \label{alg: RP_via_PP} 
\begin{algorithmic}
\small
\STATE {\textbf{Input}: $\{\lambda_{k}\} \in \mathbb{R}_{++}$ satisfying $\sum_{k=0}^{\infty} \lambda_{k}^{-1} = +\infty $ \textbf{and} $\varphi(\cdot) = g \circ \|\cdot\|_{\mathcal{V}^{\circ}}: \mathcal{X} \rightarrow \mathbb{R}$ satisfying Assumption \ref{assumption: blanket varphi}.} 
\STATE {\bf Solve for the robsut solution $x_{\mathrm{R}} = \argmin_{x\in \mathcal{X}} \varphi(x) = \Pi^{\varphi}_{\mathcal{X}}(0)$ and set $x_{0} = x_{\mathrm{R}}$.} 
\FOR{$k=0,1,...$} 
    \STATE $x_{k+1} = \argmin_{x\in \mathcal{X}} \langle a_{0},x\rangle +\lambda_{k} D_{\varphi}(x,x_{k})$ \\
\ENDFOR
\RETURN$\{x_{k}\}$ as an (approximate) robust path of $(\mathrm{RC})$. The corresponding uncertainty set radius $r_{k}$ follows the closed-form expression:
$r_{k} = \omega_{k}  \nabla g \left(\|x_{k}\|_{\mathcal{V}^{\circ}}\right),$
with $\omega_{k} = \left( \sum_{j=0}^{k-1} \lambda_{j}^{-1} \right)^{-1}.$
\end{algorithmic}
\end{algorithm}

\emph{Remark.} If the conditions of Theorem \ref{theorem: PPMP_Exact_RP} are satisfied, Algorithm \ref{alg: RP_via_PP} generates a proximal path that is an exact robust path. More generally, Algorithm 1 produces a proximal path that is an approximate robust path, where the approximation error between the two paths can be characterized via Theorem \ref{theorem: uniform bound CP and RP} and Theorem \ref{theorem: PPM as approx BPP}. We also point out that an exact proximal point step is of the same computational cost as solving a single $(\mathrm{RC})$, to lower the computational cost, the exact proximal point step of Algorithm \ref{alg: RP_via_PP} can be replaced by its computationally cheaper approximations (e.g., projected gradient descent \cite{parikh2014proximal}) to trade higher approximate errors for lower computational costs.

\section{Numerical Validations}
\label{section: application}
In the previous sections, we developed a theory for robust paths and related optimization paths; the theoretical insights are operationalized as Algorithm \ref{alg: RP_via_PP}. In this section, we validate the results of our theorems and the performance of Algorithm \ref{alg: RP_via_PP} via numerical experiments. 

The first experiment in Section \ref{sec: exp-portfolio} on portfolio optimization illustrates how Algorithm \ref{alg: RP_via_PP} generates an entire set of approximate, sometimes exact, Pareto efficient portfolios via a single {\ProximalPathFull}. We compare Algorithm \ref{alg: RP_via_PP} generated portfolios with exact Pareto efficient portfolios under three types of setups: a hyperplane feasible region, a polyhedrally monotone feasible region, and a more general feasible region with extra practical investment constraints. 
The hyperplane feasible region shows as special cases of Theorem \ref{theorem: PPMP_Exact_RP}, under affine subspace feasible regions or unconstrained problems, the {\ProximalPathFull}s are exact robust paths. 
The polyhedrally monotone feasible region confirms under the conditions of Theorem \ref{theorem: PPMP_Exact_RP}, the {\ProximalPathFull}s are exact robust paths. The general feasible region setup confirms the results of Theorem \ref{theorem: uniform bound CP and RP} and Theorem \ref{theorem: PPM as approx BPP} that {\ProximalPathFull}s in general can be good approximate robust paths.
As a side product, our results extend the classical Two-Fund Theorem in finance \citep{markowitz2008portfolio}, with our Algorithm \ref{alg: RP_via_PP} accommodating practical trading constraints while the original theorem in finance does not.

The second experiment in Section \ref{sec: exp_ADML} explores settings that significantly deviate from our theoretical assumptions, with highly nonlinear loss functions in deep learning. Even here, our theory and algorithm (correctly initializing, and sometimes restarting standard gradient-based methods) lead to orders of magnitude computational speedup while maintaining strong worst-case and nominal prediction performances. 

\subsection{First Experiment: Portfolio Optimization}
\label{sec: exp-portfolio}

\begin{figure}[t]
    \centering
    \begin{subfigure}[t]{0.33\linewidth}
        \centering
        \includegraphics[width=\linewidth]{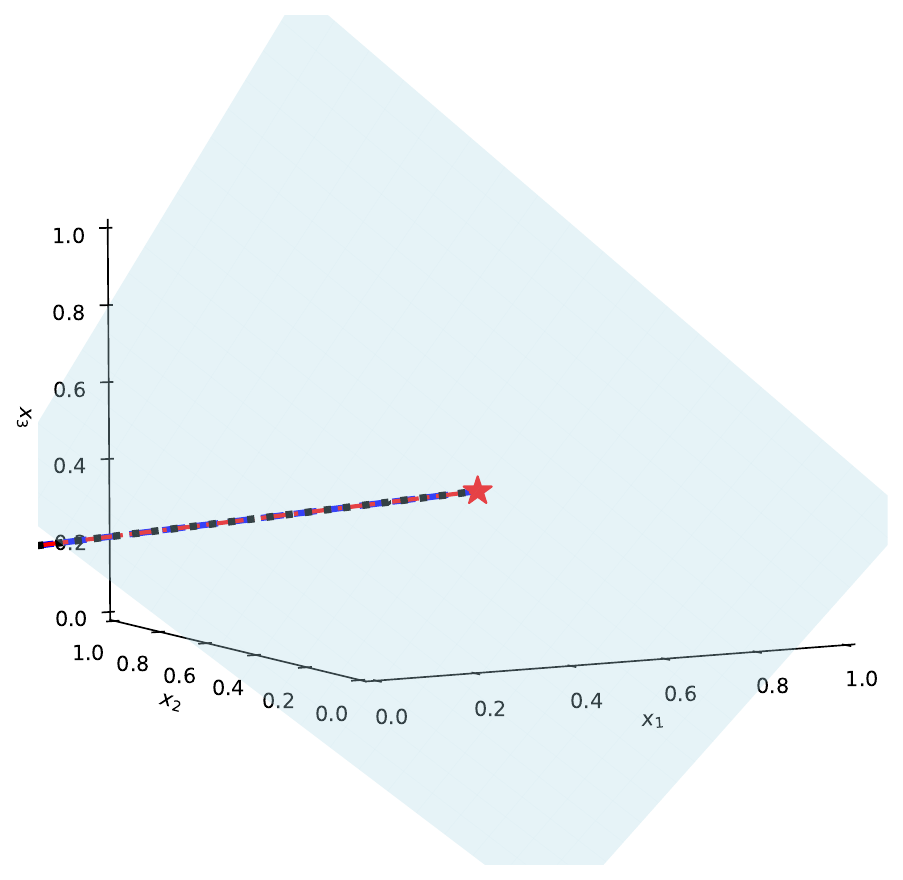}
        \caption*{(a1) Hyperplane: $\mathcal{X} = \{x\in \mathbb{R}^{3}: \langle \mathbf{1},x\rangle=1 \}$.}
    \end{subfigure}%
    \begin{subfigure}[t]{0.33\linewidth}
        \centering
        \includegraphics[width=\linewidth]{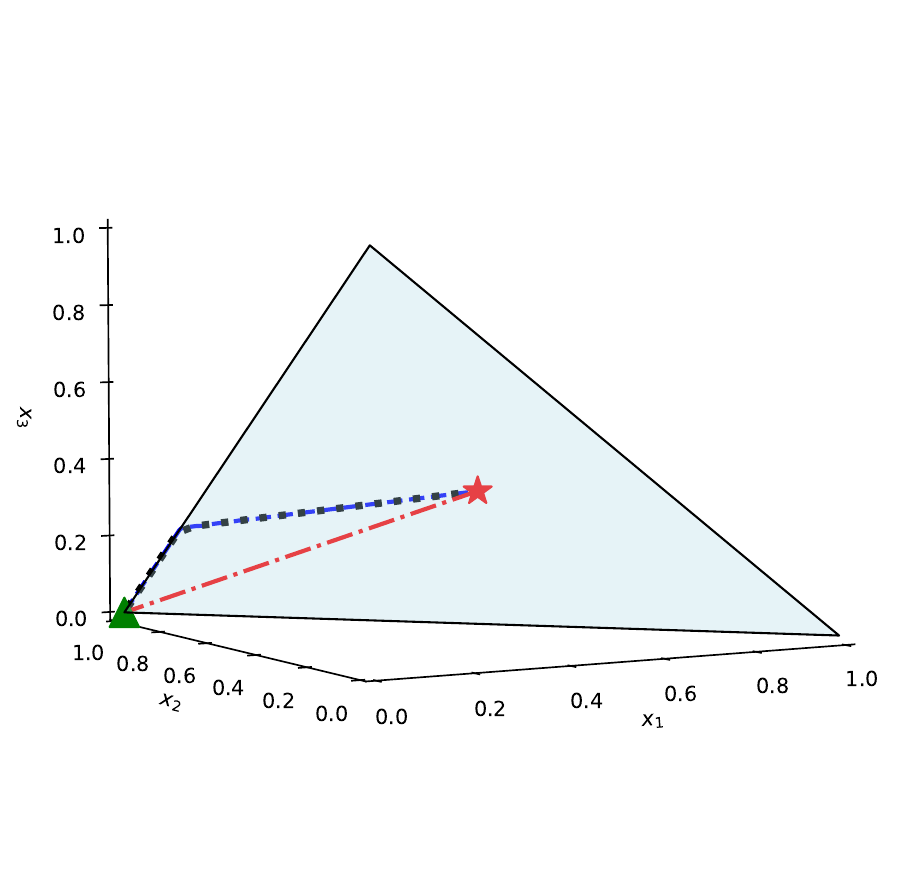}
        \caption*{(a2) Simplex: $\mathcal{X} = \{x\in \mathbb{R}^{3}: \langle \mathbf{1},x\rangle=1, \ x\geq 0\}$.}
    \end{subfigure}
    \begin{subfigure}[t]{0.33\linewidth}
        \centering
        \includegraphics[width=\linewidth]{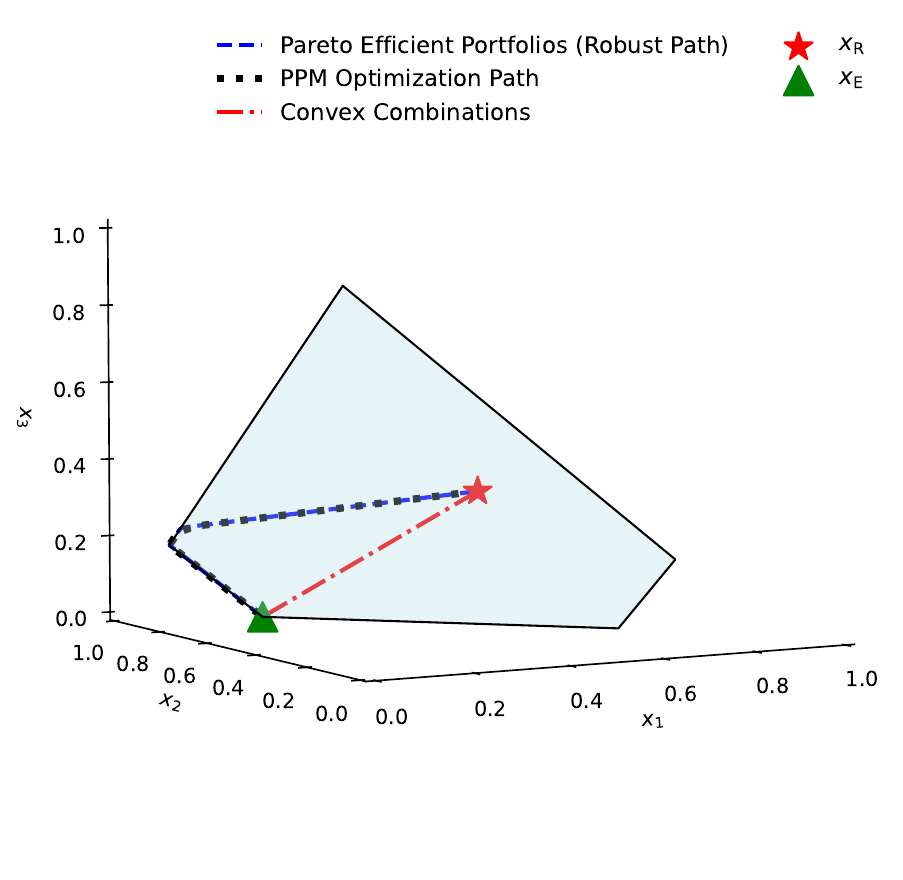}
        \caption*{ (a3) General:  $\mathcal{X} = \{x\in \mathbb{R}^{3}: \langle \mathbf{1},x\rangle=1, \  x\geq [0.0,0.1,0.0], \  x\leq [0.7,0.8,1.0]\}$.}
    \end{subfigure}
        \begin{subfigure}[t]{0.33\linewidth}
        \centering
        \includegraphics[width=\linewidth]{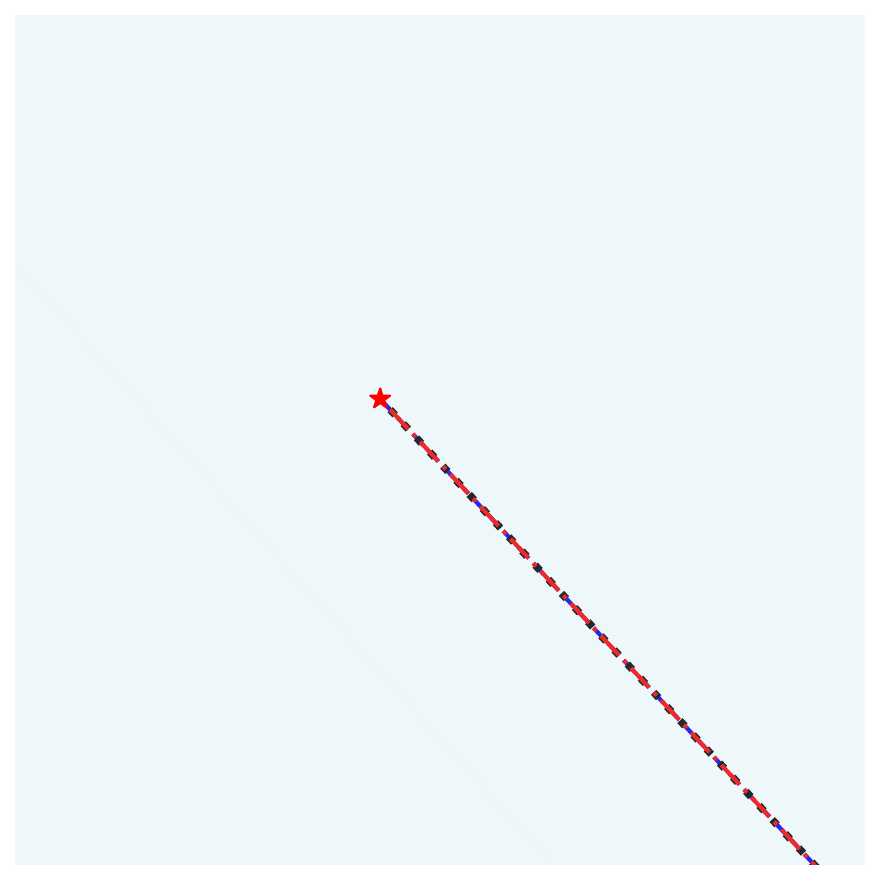}
        \caption*{(b1) Hyperplane: $\mathcal{X} = \{x\in \mathbb{R}^{4}: \langle \mathbf{1},x\rangle=1 \}$.}
    \end{subfigure}%
    \begin{subfigure}[t]{0.33\linewidth}
        \centering
        \includegraphics[width=\linewidth]{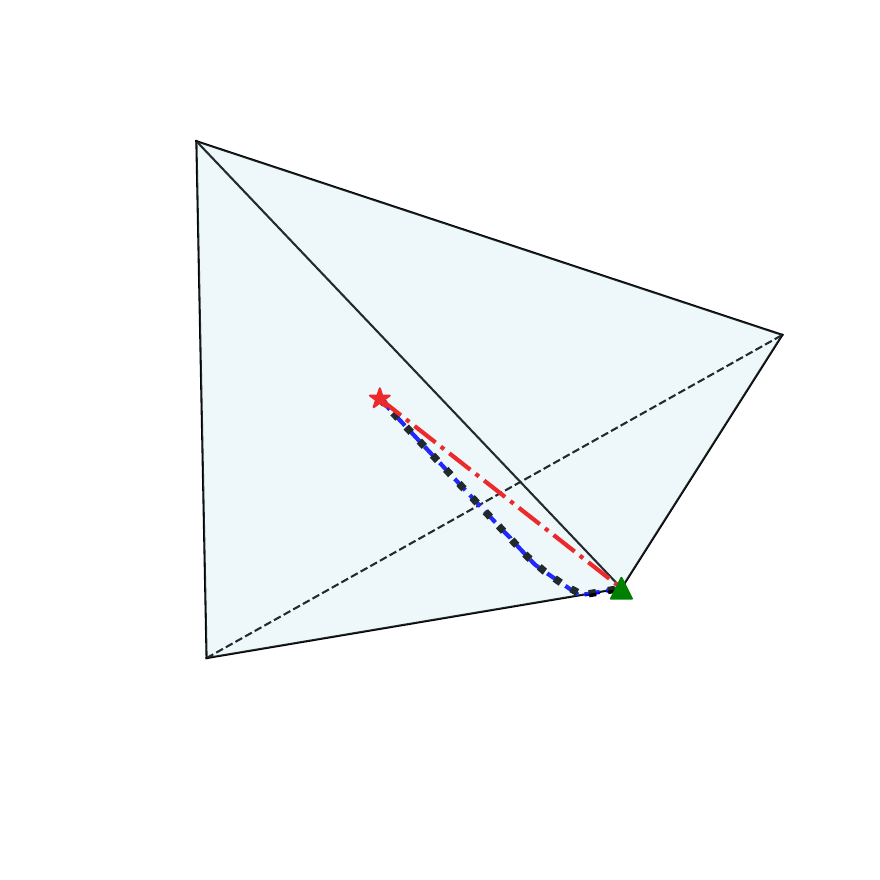}
        \caption*{(b2) Simplex:  $\mathcal{X} = \{x\in \mathbb{R}^{4}: \langle \mathbf{1},x\rangle=1, \ x\geq 0\}$.}
    \end{subfigure}
    \begin{subfigure}[t]{0.33\linewidth}
        \centering
        \includegraphics[width=\linewidth]{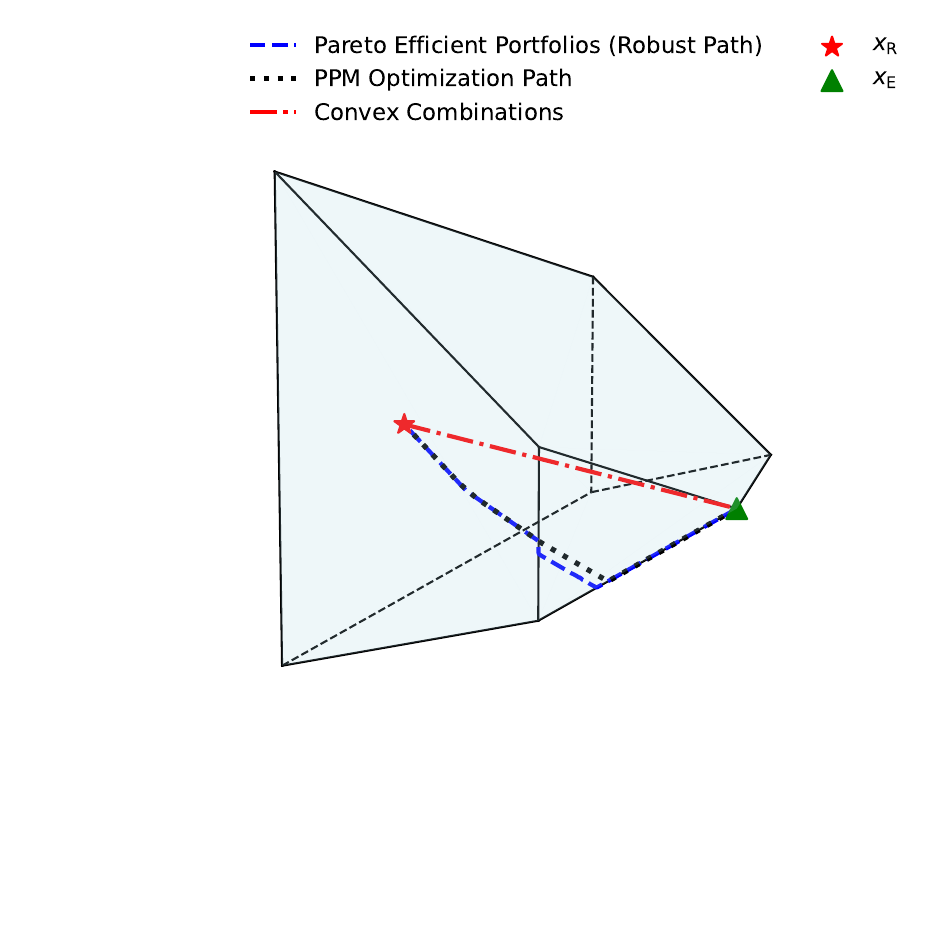}
        \caption*{(b3) General:  $\mathcal{X} = \{x\in \mathbb{R}^{4}: \langle \mathbf{1},x\rangle=1. \ x\geq [0.0,0.1,0.0,0.0] ,\ x \leq [0.6,0.6,1.0,1.0]\}$.}
    \end{subfigure}
    \caption{{\ProximalPathFullCap}s and convex combinations in the solution space, as approximate Pareto efficient portfolios (equivalently, robust paths) of problem (\ref{eq: robust_portfolio}): under hyperplane, simplex, and general feasible regions. The left two columns represent cases where our theoretical results predict a precise alignment between the {\ProximalPathFull}s and robust paths (black and blue lines). Classical Two-Fund Theorem (red line) works for the leftmost column (unconstrained case), but not for constrained cases. In the rightmost column, our Theorem \ref{theorem: PPM as approx BPP} predicts a small gap between the robust paths (blue line) and {\ProximalPathFull}s (black line) since the feasible region is no longer polyhedrally monotone. We observe exactly that, and also in addition a degeneracy where the top right case (a3) shows blue and black lines actually coincide. On the contrary, the classical Two-Fund Theorem (red line) can not generate useful portfolios to match Pareto efficiency.}
    \label{fig: two-fund generalized}
\end{figure}

In portfolio optimization, we are concerned with constructing a portfolio from $n$ risky assets. The return of the $n$ assets is modeled by a random vector, $\alpha$. We assume that from historical data, we can estimate the expectation and the covariance matrix of $\alpha$ to be $\mu$ and $\Sigma$. It is known that the classical Markowitz mean-variance portfolio optimization problem can be cast in the form of our $(\mathrm{RC})$ with an ellipsoidal uncertainty set $\mathcal{U}(r) = \{\mu + \xi \in \mathbb{R}^{n}: \ \|\Sigma^{-1/2} \xi \|_{2} \leq r \}$ \citep{natarajan2009constructing}:
\begin{equation}   
\label{eq: robust_portfolio}
\min_{x \in {\mathcal{X}}} \max_{\alpha \in \mathcal{U}(r)} - \langle \alpha,x \rangle 
\end{equation}

\subsubsection{{\ProximalPathFullCap}s are Almost Pareto Efficient Portfolios: Visualization in Solution Space.} We compute the exact Pareto efficient portfolios (equivalently, the robust path) of problem (\ref{eq: robust_portfolio}). In addition, we compute two approximations of the set of Pareto efficient portfolios: (i) Two-Fund Theorem \citep{markowitz2008portfolio}: convex combinations of the min-variance $x_{\mathrm{R}}$ and the max-return portfolio, $x_{\mathrm{E}}$ and (ii) Algorithm \ref{alg: RP_via_PP}: {\ProximalPathFull} initialized at $x_{\mathrm{R}}$ and converging to $x_{\mathrm{E}}$. We run the above experiment under increasingly more general asset weights feasible regions, $\mathcal{X}$: hyperplane, simplex, and simplex with additional trading constraints. 

Figure \ref{fig: two-fund generalized} presents the results on two small instances where a portfolio is constructed with three (top row) and four (bottom row) risky assets. The results give the following insights: 
\begin{itemize}
    \item[(i)] Under general feasible regions $\mathcal{X}=\{x\in \mathbb{R}^{n}: \langle 1, x \rangle=1, x_{\mathrm{lb}}\leq x\leq x_{\mathrm{ub}}\}$: As depicted in Figure \ref{fig: two-fund generalized}, (a3) and (b3), the {\ProximalPathFull}s are higher quality approximations of Pareto efficient portfolios than convex combinations. For instance, Figure \ref{fig: two-fund generalized}, (a3) shows a degenerate instance where the {\ProximalPathFull} is not monotone on $\mathcal{X}$ (entered one edge from another edge), but remains an exact set of Pareto efficient portfolios. Figure \ref{fig: two-fund generalized}, (b3) shows a general instance where the {\ProximalPathFull} is a set of approximate Pareto efficient portfolios, where the approximation error bound can be characterized by Theorem \ref{theorem: uniform bound CP and RP} and Theorem \ref{theorem: PPM as approx BPP}.
    \item[(ii)] Under simplex feasible regions $\mathcal{X}=\{x\in \mathbb{R}^{n}: \langle 1, x \rangle=1, x\geq 0\}$: As shown in Figure \ref{fig: two-fund generalized}, (a2) and (b2), the set of Pareto efficient portfolios has a piecewise linear structure, under which the Two-Fund Theorem no longer holds and convex combinations only generate poor approximations of the Pareto efficient portfolios. In contrast, by Proposition \ref{prop: monotone_BPP_is_PPM}, the {\ProximalPathFull}s are monotone on simplices; hence, the {\ProximalPathFull}s are exact Pareto efficient portfolios.
    \item[(iii)] Under hyperplane feasible regions $\mathcal{X}=\{x\in \mathbb{R}^{n}: \langle 1, x \rangle=1\}$: the Two-Fund Theorem \citep{markowitz2008portfolio} states that the set of Pareto efficient portfolios can be constructed exactly as convex combinations of the min-variance portfolio, $x_{\mathrm{R}}$, and the max-return portfolio, $x_{\mathrm{R}}$. Our Theorem \ref{theorem: RP_CP_PP are Bregman Projections} gives a geometric proof of the Two-Fund Theorem, as depicted in Figure \ref{fig: two-fund generalized}, (a1) and (b1): under hyperplane feasible regions, the set of Pareto efficient portfolios is a line segment in $\mathbb{R}^{n}$ which can be generated as convex combinations of any two Pareto efficient portfolio including $x_{\mathrm{R}}$ and $x_{\mathrm{E}}$. In addition, by Proposition \ref{prop: monotone_BPP_is_PPM}, the {\ProximalPathFull}s are monotone on hyperplanes; hence, the {\ProximalPathFull}s are exact Pareto efficient portfolios.

\end{itemize}

\subsubsection{Performance of {\ProximalPathFullCap} Portfolios: Visualization in the Objective Space.} Next, we test the performance of the {\ProximalPathFull} generated portfolios against the exact Pareto efficient portfolios, as measured by the worst-case return and nominal case return. 

\textbf{Experiment setup.} We construct approximate Pareto-efficient portfolios with {\ProximalPathFull}s. We consider two general feasible regions $\Delta^{n} = \{x\in \mathbb{R}^{n}: \langle \mathbf{1},x\rangle=1, x\geq 0\}$ and the feasible region of Markowitz++ model as introduced by \cite{boyd2024markowitzportfolioconstructionseventy}.  We first construct portfolios with in-sample historical stock return data ($20$ stocks within S$\&$P 500 from 2021-01-01 to 2023-12-30), before testing the nominal and worst-case returns of each portfolio on out-of-sample stock return data (same stocks from 2024-01-01 to 2024-08-01).  As the benchmark, we first construct exact Pareto efficient portfolios by solving problem (\ref{eq: robust_portfolio}) exactly multiple times each under a different $r$. Then we run a single pass of {\ProximalPointMethodAbb} to generate approximate Pareto efficient robust portfolios, \ie, we use the minimum-variance (most robust) portfolio with $r=\infty$ to initialize the {\ProximalPathFull} for solving the nominal portfolio optimization problem with $r=0$, the {\ProximalPathFull} generates approximate Pareto-efficient portfolios. 

\begin{figure}[ht]
    \centering
    \includegraphics[width=\linewidth]{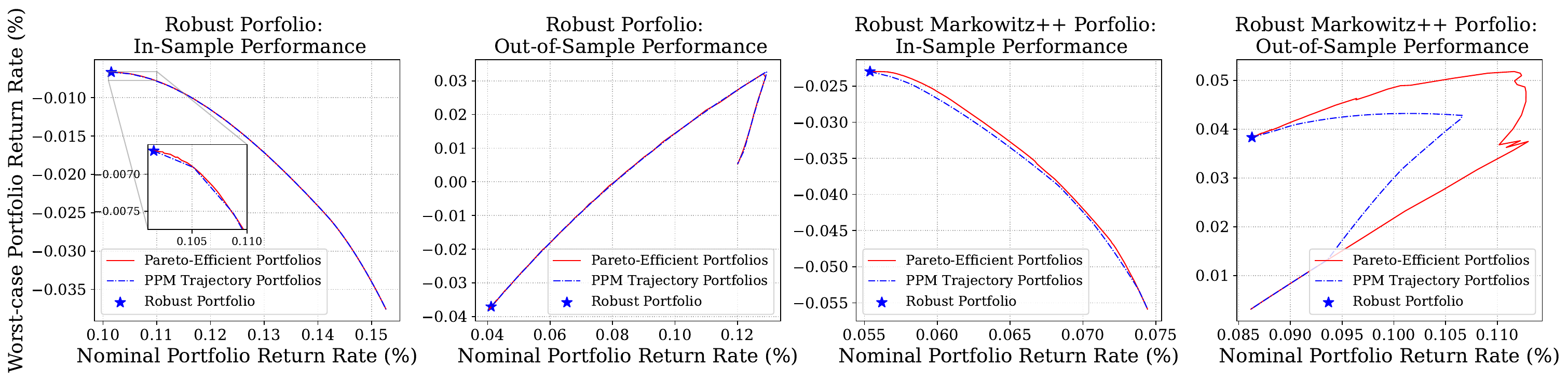}
    \caption{In-sample and out-of-sample performances (as measured by the nominal and worst-case returns) of the exact Pareto-efficient portfolios v.s. {\ProximalPathFull} generated approximate Pareto-efficient path portfolios.}
    \label{fig:Portfolio Pareto-frontier}
\end{figure}

\textbf{Results.} The experiment results are presented in Figure \ref{fig:Portfolio Pareto-frontier}, the {\ProximalPathFull} generated approximate Pareto-efficient portfolios match closely the exact Pareto-efficient portfolios in performance as measured by the nominal and the worst-case return rates.

\subsection{Second Experiment: Obtaining Multiple Robust Solutions in Deep Learning}
\label{sec: exp_ADML}

\noindent \textbf{Adversarial training as robust optimization.} The goal in adversarially robust deep learning is to learn networks that are robust against adversarial attacks (\ie, perturbations on the input examples that aim to deteriorate the accuracy of classifiers). A common strategy to robustify networks is adversarial training, which can be formulated as the following robust optimization problem \citep{madry2018towards},
\begin{equation}
\label{eq: adv_training}
    \min_{\theta} \mathbb{E}_{(x,y) \sim \mathcal{D}} \left[\max_{\xi \in \Xi(r,\mathcal{V})}\ell(f_{\theta}(x + \xi), y) \right],
\end{equation}
where $\mathcal{D}$ is the distribution generating pairs of examples $x\in\mathbb{R}^{d}$ and labels $y \in [c]$, $f_{\theta}$ is a neural network parameterized by $\theta$, $\xi$ is the perturbation/attack on the input data, within a perturbation set $\Xi(r,\mathcal{V})$, and $\ell$ is the lost function.  Standard adversarial training methods \citep{madry2018towards, Wong2020Fast} approximately solve problem (\ref{eq: adv_training}). If adversarial training under a fixed $r$ can be computed in $O(T_{\mathrm{at}})$ time, generating $n$ adversarially robust models under varying $r$ (\ie, a robust path) takes $O(nT_{\mathrm{at}})$.

As an adaptation of Algorithm \ref{alg: RP_via_PP} for problem (\ref{eq: adv_training}), we propose Algorithm \ref{alg: PERO via FOM} to compute an approximate robust path of problem (\ref{eq: adv_training}) in only two algorithmic passes: a single \emph{adversarial training} followed by a single \emph{standard training} with approximate {\ProximalPointMethodAbb} initialized with $x_{\mathrm{R}}$; finally the approximate {\ProximalPointMethodAbb} iterates are an approximate robust path of problem (\ref{eq: adv_training}). Assume one step of an approximate {\ProximalPointMethodAbb} step in standard training costs $O(T_{\mathrm{ppm}})$, generating $n$ adversarially robust models under varying $r$ via algorithm \ref{alg: PERO via FOM} takes $O(T_{\mathrm{at}} + nT_{\mathrm{ppm}})$.

\begin{algorithm}[ht!]
\caption{Approximate Robust Path of Problem (\ref{eq: adv_training}) via {\ProximalPointMethodFullCap}} \label{alg: PERO via FOM}
\begin{algorithmic}
\small
\STATE {\textbf{Input}: $\{\lambda_{k}\} \in \mathbb{R}_{++}$ satisfying $\sum_{k=0}^{\infty} \lambda_{k}^{-1} = +\infty $ \textbf{and} $\varphi: \mathcal{X} \rightarrow \mathbb{R}$ satisfying Assumption \ref{assumption: blanket varphi}.} 
\STATE {\bf Solve problem (\ref{eq: adv_training}) under a large radius $r=\overline{r}$ for $\theta_{\mathrm{R}}$ and set $\theta_{0} = \theta_{\mathrm{R}}$.} 
\FOR{$k=0,1,...$} 
    \STATE $\theta_{k+1} \text{{$\approx$}} \argmin_{\theta} \ell(f_{\theta}(x), y) +\lambda_{k} D_{\varphi}(\theta,\theta_{k})$ \\
\ENDFOR
\RETURN$\{\theta_{k}\}$ as an approximate robust path of problem (\ref{eq: adv_training}). 
\end{algorithmic}
\end{algorithm}

\begin{figure}[t]
    \centering
    \includegraphics[width=\linewidth]{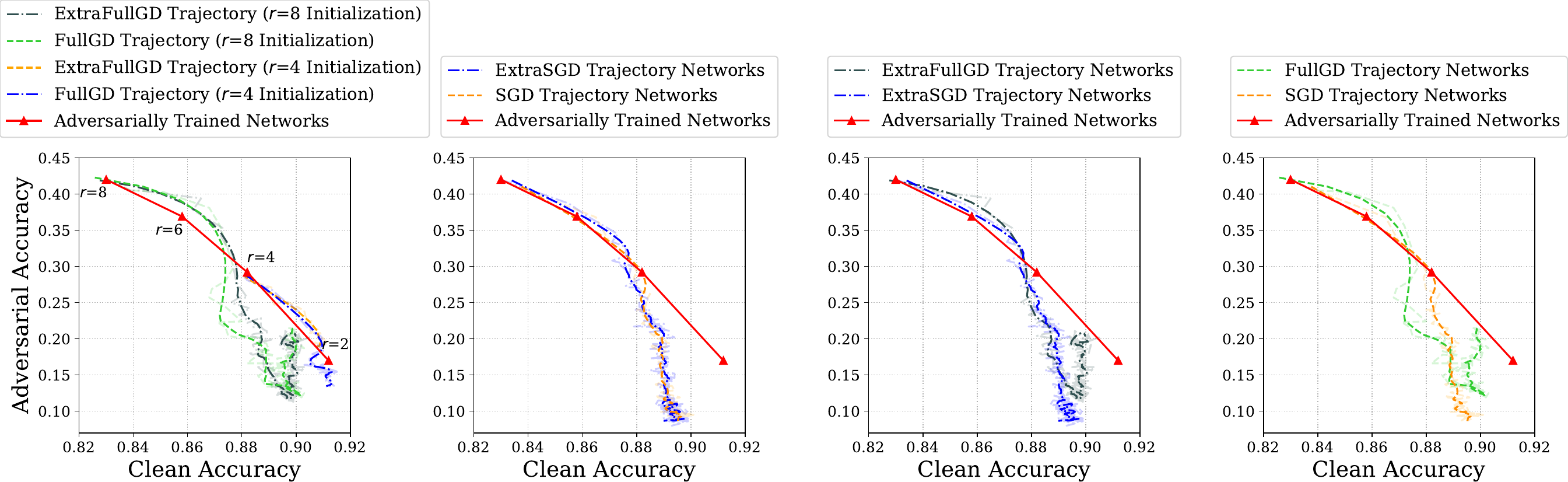}
    \caption{Robust Path of Problem (\ref{eq: adv_training}): Algorithm 1 v.s. FGSM. Algorithm 1 is equipped with the Extra-gradient descent (ExtraFullGD) and the vanilla gradient descent (FullGD) methods as approximate PPM; initialization at $\overline{r}=8$ or $4$.}
    \label{fig:Adv_ML Pareto-frontier}
\end{figure}

\begin{table}[t]
\centering
\small
\begin{tabular}{p{4cm} p{4.5cm} p{6cm}}
\toprule
Method & Time per robust model (min) & Time to solve $N=100$ robust models (hrs)\\
\midrule
Our Method: Algorithm \ref{alg: PERO via FOM} & $0.25$ ($T_{\mathrm{ppm}}$)   & $0.64$ ($T_{\mathrm{at}} + nT_{\mathrm{ppm}}$)   \\
FGSM \citep{Wong2020Fast} & $15.12$ ($T_{\mathrm{at}}$) & $ 25.2$ ($nT_{\mathrm{at}}$) \\
\bottomrule
\end{tabular}
\caption{Computation Cost: Algorithm \ref{alg: PERO via FOM} v.s. FGSM}
\label{table:Adv_ML_Cost_Compare}
\end{table}

\noindent \textbf{Experiment result.} We set up problem (\ref{eq: adv_training}) as training robust image classifiers with the CIFAR10 dataset and a PreAct ResNet18 architecture. As the benchmark for algorithm \ref{alg: PERO via FOM}, we first adversarially train networks with the state-of-the-art FGSM \citep{Wong2020Fast} under different perturbation set radii, $r$. Then we run Algorithm \ref{alg: PERO via FOM} to generate another approximate robust path of problem (\ref{eq: adv_training}). The performance (clean accuracy v.s. adversarial accuracy) of Algorithm \ref{alg: PERO via FOM} computed robust paths is comparable to that of FGSM (Figure \ref{fig:Adv_ML Pareto-frontier}). At the same time, the computation cost of our Algorithm \ref{alg: PERO via FOM} is significantly lower than that of FGSM (Table \ref{table:Adv_ML_Cost_Compare}).

\section{Conclusion}

We discuss three pieces of future research that are out of the scope of this paper, but are immediate and important directions for possible future papers:
\begin{itemize}
   \item \emph{Nonlinear objective functions.} 
   Although our theory is on robust optimization problems with linear objectives as per its closest robust optimization literature of \cite{Bertsimas2004, Iancu2014, Stability_2}, our results can be generalized to convex robust optimization problems via an epigraph reformulation, where the geometrical view of Theorem \ref{theorem: RP_CP_PP are Bregman Projections} and all subsequent results can be restated through epigraphical projections. In addition, as demonstrated in Section \ref{sec: exp_ADML}, the computational speedup of our framework is significant in recovering the robust path of nonlinear (nonconvex) robust optimization problems, where the exact proximal point method is replaced by its computationally cheaper approximations. This also leads to the next research direction.
   \item \emph{Approximate {\ProximalPointMethodFull}: trade-off approximation quality with computation cost.} In the practical use of Algorithm 1, especially when applied to nonlinear robust optimization problems, the exact proximal point steps can be replaced by its computationally cheaper approximations (\eg, projected gradient descent \citep{parikh2014proximal}, extra-gradient, optimistic gradient \citep{mokhtari2020unified}) to trade a higher approximation error to the robust path for a lower computational cost. The approximation error between the proximal path and its cheaper algorithmic approximations can be combined directly with our approximation error between the robust path and the proximal path.

    \item \emph{Automated tuning of robustness and efficiency trade-off.} Our Algorithm 1 recovers the entire robust path. What remains to be answered is which robust solution from the robust path should be deployed? On a high level, our results take the crucial step towards addressing this problem by reducing the high-dimensional solution search space to the one-dimensional robust path recoverable via a proximal path. A natural follow-up work is the algorithmic automation for selecting the robust solution with an appropriate robustness and efficiency trade-off on the now reduced search space of the robust path attainable via a proximal path.
\end{itemize}

\bibliographystyle{informs2014} 
\bibliography{bibliography} 

\newpage

\begin{APPENDICES}

\section{Proof of Lemma \ref{lemma: Breg_Proj_calc_1}}
\label{APP: proof of Breg_Proj_calc_1}
Denote $x^{+} = \Pi^{\varphi}_{\mathrm{Aff}(\mathcal{C})}(x)$ and $x^{++} = \Pi^{\varphi}_{\mathcal{C}} \circ \Pi^{\varphi}_{\mathrm{Aff}(\mathcal{C})}(x)$. By Lemma \ref{lemma: Bregman Projection VI},
\begin{subequations}
\label{eq: def x^plus}
\begin{align}
\langle \nabla \varphi(x) - \nabla \varphi(x^{+}), z-x^{+}  \rangle &= 0 , \quad \forall z \in \mathcal{\mathrm{Aff}(C)}, \\
\Rightarrow \ \langle \nabla \varphi(x) - \nabla \varphi(x^{+}), z-x^{+}  \rangle &= 0 , \quad \forall z \in \mathcal{C},
\end{align}
\end{subequations} 
and
\begin{equation}
\label{eq: def x^plusplus}
\langle \nabla \varphi(x^{+}) - \nabla \varphi(x^{++}), z-x^{++}  \rangle \leq 0 , \quad \forall z \in \mathcal{C}.
\end{equation}
Given $x^{++} \in \mathcal{C}$ and (\ref{eq: def x^plus}) we also have
\begin{equation}
\label{eq: def x^plus x^plusplus in C}
    \langle \nabla \varphi(x) - \nabla \varphi(x^{+}), x^{++}-x^{+}  \rangle = 0.
\end{equation}
Adding (\ref{eq: def x^plus}) and (\ref{eq: def x^plusplus}) we have for all $z\in \mathcal{C}$,
$$
\begin{aligned}
    \langle \nabla \varphi(x), z\rangle + \left( -  \langle \nabla \varphi(x), x^{+} \rangle -  \langle \nabla \varphi(x^{+}), x^{++} \rangle + \langle \nabla \varphi(x^{+}), x^{+} \rangle \right) \\
    - \langle \nabla \varphi(x^{++}), z \rangle + \langle \nabla \varphi(x^{++}), x^{++} \rangle &\leq 0 \\
    \langle \nabla \varphi(x), z\rangle - \langle \nabla \varphi(x), x^{++} \rangle - \langle \nabla \varphi(x^{++}), z \rangle + \langle \nabla \varphi(x^{++}), x^{++} \rangle &\leq 0 \\
    \langle \nabla \varphi(x) - \nabla \varphi(x^{++}), z- x^{++} \rangle & \leq 0,
\end{aligned}
$$
where the second inequality is due to (\ref{eq: def x^plus x^plusplus in C}). Finally by Lemma \ref{lemma: Bregman Projection VI}, we have $x^{++} = \Pi^{\varphi}_{\mathcal{C}}(x)$.  \hfill \Halmos

\section{Proofs of Lemmas \ref{lemma: dual_RC} and \ref{lemma: unique robust path}}
\label{sec: proof proposition: unique robust path}
\begin{proof}{Proof of Lemma \ref{lemma: dual_RC}.}
Given $\Xi = \Xi(r,\mathcal{V})=  \left\{\xi \in \mathbb{R}^{n}:    \|\xi\|_{\mathcal{V}} \leq r  \right\},$
$$
\begin{aligned}
    &\min_{x \in \mathcal{X}} \max_{a \in \mathcal{U}} \langle a,x \rangle \\
    = & \min_{x \in \mathcal{X}} \langle a_{0},x \rangle   + \max_{\xi \in \mathbb{R}^{n}} \left\{\langle \xi, x\rangle : \|\xi\|_{\mathcal{V}} \leq r  \right\}  \\
    = & \min_{x \in \mathcal{X}} \langle a_{0},x \rangle   + \max_{\phi \in \mathbb{R}^{n}} \left\{ r \langle \phi, x\rangle : \phi \in \mathcal{V}\right\}  \\
    = &  \min_{x \in \mathcal{X}} \langle a_{0},x \rangle   +r \|x\|_{\mathcal{V}^{\circ}}.
\end{aligned}
$$ 
The last equality is due to Assumption \ref{assumption: blanket V}: $0\in \mathrm{int}(\mathcal{V})$ and \cite{Rockafellar_Convex_Analysis} Theorem 14.5.
\hfill \Halmos
\end{proof}

\begin{proof}{Proof of Lemma \ref{lemma: unique robust path}.}
For any $\omega > 0$, by the definition of $x_{\mathrm{R}}'(\omega, \mathcal{V})$, 
$$\langle  a_0 + \omega \nabla g (\|x_{\mathrm{R}}'\left(\omega, \mathcal{V})\|_{\mathcal{V}^{\circ}} \right) \nabla\|x_{\mathrm{R}}'\left(\omega, \mathcal{V} \right)\|_{\mathcal{V}^{\circ}}, x - x_{\mathrm{R}}'\left(\omega, \mathcal{V} \right) \rangle \geq 0, $$
given $r(\omega) = \omega  \nabla g (\|x_{\mathrm{R}}'\left(\omega, \mathcal{V})\|_{\mathcal{V}^{\circ}} \right)$, we have
$$
\begin{aligned}
 & \ \langle  a_0 + r(\omega)  \nabla\|x_{\mathrm{R}}'\left(\omega, \mathcal{V} \right)\|_{\mathcal{V}^{\circ}}, x - x_{\mathrm{R}}'\left(\omega, \mathcal{V} \right) \rangle \geq 0 \\
 \Leftrightarrow & \  x_{\mathrm{R}}'(\omega, \mathcal{V}) \in \argmin_{x\in \mathcal{X}} \langle a_{0}, x \rangle + r(\omega)\|x\|_{\mathcal{V}^{\circ}}. 
\end{aligned}
$$
By Lemma \ref{lemma: dual_RC},
$$x_{\mathrm{R}}'(\omega, \mathcal{V}) \in \argmin_{x \in \mathcal{X}}  \max_{\xi  \in\Xi(r(\omega),\mathcal{V})} \langle a_{0} +  \xi,x \rangle. $$ \hfill \Halmos
\end{proof}

\section{$\kappa$-Expansiveness Example}
\label{sec:kappa-expansiveness-example}

When $\varphi = \frac{1}{2}\|x\|^{2}_{2}$, the usual Euclidean projection $\Pi^{\varphi}_{\mathcal{S}}$ is 1-expansive. More generally, for smooth and strongly convex $\varphi$, it is easy to verify the following result.
\begin{definition}
    $\varphi$ is $L\text{-smooth}$ and $\mu\text{-strongly convex}$ w.r.t. some norm $\|\cdot\|$ if 
    $$\frac{\mu}{2} \|x - y\|^{2}\leq D_{\varphi}(x,y) \leq \frac{L}{2} \|x - y\|^{2}, \quad \forall x,y \in \mathrm{dom} (\varphi).$$
\end{definition}
\begin{proposition}
\label{proposition: bregman proj expansiveness 1}
 Assume $\varphi$ is $L\text{-smooth}$ and $\mu\text{-strongly convex}$ w.r.t. some norm $\|\cdot\|$, then the induce Bregman projection $\Pi^{\varphi}_{\mathcal{S}}$ is $\kappa$-expansive with $\kappa = \left( \frac{L}{\mu} \right)^{3}$.
\end{proposition}
\begin{proof}{Proof.}
We begin by proving
   $D_{\varphi}\left(\Pi^{\varphi}_{\mathcal{S}}(x), \Pi^{\varphi}_{\mathcal{S}}(y) \right) \leq  (L/\mu)^{3} \cdot D_{\varphi} \left(x,y \right), \ \forall x,y \in\mathrm{int}(\mathrm{dom}(\varphi))$:

Denote $x^{+} = \Pi^{\varphi}_{\mathcal{X}}(x)$ and $y^{+} = \Pi^{\varphi}_{\mathcal{X}}(y)$, by Lemma \ref{lemma: Bregman Projection VI}, we have
\begin{equation}
\langle \nabla \varphi(x) - \nabla \varphi(x^{+}), z-x^{+}  \rangle \leq 0 , \quad \forall z \in \mathcal{X},
\end{equation}
and
\begin{equation}
\langle \nabla \varphi(y) - \nabla \varphi(y^{+}), z-y^{+}  \rangle \leq 0 , \quad \forall z \in \mathcal{X}.
\end{equation}
Given $x^{+} \in \mathcal{X}$ and $y^{+} \in \mathcal{X}$, 
\begin{equation}
\label{eq:lemma9_3}
\langle \nabla \varphi(x) - \nabla \varphi(x^{+}), y^{+}-x^{+}  \rangle \leq 0 ,
\end{equation}
and
\begin{equation}
\label{eq:lemma9_4}
\langle \nabla \varphi(y) - \nabla \varphi(y^{+}), x^{+}-y^{+}  \rangle \leq 0 .
\end{equation}
Combining ($\ref{eq:lemma9_3}$) and ($\ref{eq:lemma9_4}$),
$$
    \langle \nabla \varphi(x^{+}) - \nabla \varphi(y^{+}), x^{+}-y^{+}  \rangle \leq    \langle \nabla \varphi(x) - \nabla \varphi(y), x^{+}-y^{+}  \rangle.
$$
By the strong convexity of $\varphi$,
\begin{subequations}
\begin{align}
    \mu \|x^{+} - y^{+}\|^{2} &\leq   \langle \nabla \varphi(x) - \nabla \varphi(y), x^{+}-y^{+}  \rangle \\ 
    \mu \|x^{+} - y^{+}\|^{2} &\leq   \| \nabla \varphi(x) - \nabla \varphi(y) \|_{*} \cdot \|x^{+}-y^{+}\|   \\
    \|x^{+} - y^{+}\|^{2} & \leq \mu^{-2}  \| \nabla \varphi(x) - \nabla \varphi(y) \|_{*}^{2}.
\end{align}
\end{subequations}
By the smoothness of $\varphi$,
\begin{equation}
\label{eq:lemma9_euclidean_proj_expansiveness}
    \|x^{+} - y^{+}\|^{2} \leq (\mu/L)^{-2} \|x-y \|^{2}.
\end{equation}
By the definition of smoothness and strong convexity,
$$
 D_{\varphi}(x^{+}, y^{+}) \leq \frac{L}{2}\|x^{+} - y^{+}\|^{2} \text{\quad and \quad }  \frac{\mu}{2} \|x - y\|^{2}\leq D_{\varphi}(x,y),
$$
Together with (\ref{eq:lemma9_euclidean_proj_expansiveness}), 
\begin{equation}
    D_{\varphi}(x^{+}, y^{+}) \leq (L/\mu)^{3}  D_{\varphi}(x, y).
\end{equation}

Next, we prove that denote $d = y-x$, $D_{\varphi}\left(\Pi^{\varphi}_{\mathcal{S}+d}(x), \Pi^{\varphi}_{\mathcal{S}}(x) \right) \leq  (L/\mu)^{3} \cdot D_{\varphi} \left(x,y \right), \ \forall x,y \in\mathrm{int}(\mathrm{dom}(\varphi))$: 
Denote $x^{++} = \Pi^{\varphi}_{\mathcal{X}+d}(x)$ and $x^{+} = \Pi^{\varphi}_{\mathcal{X}}(x)$, by lemma \ref{lemma: Bregman Projection VI}, 
\begin{equation}
\langle \nabla \varphi(x) - \nabla \varphi(x^{++}), z-x^{++}  \rangle \leq 0 , \quad \forall z \in \mathcal{X}+d,
\end{equation}
and
\begin{equation}
\langle \nabla \varphi(x) - \nabla \varphi(x^{+}), z-x^{+}  \rangle \leq 0 , \quad \forall z \in \mathcal{X}.
\end{equation}
Given $x^{+} \in \mathcal{X}$ and $x^{++} \in \mathcal{X}+d$, or equivalently $x^{+}+d \in \mathcal{X}+d$ and $x^{++} -d \in \mathcal{X}$, we have
\begin{equation}
\label{eq:lemma10_3}
\langle \nabla \varphi(x) - \nabla \varphi(x^{++}), x^{+}+d-x^{++}  \rangle \leq 0 ,
\end{equation}
and
\begin{equation}
\label{eq:lemma10_4}
\langle \nabla \varphi(x) - \nabla \varphi(x^{+}), x^{++} - d -x^{+}  \rangle \leq 0 .
\end{equation}
Combining ($\ref{eq:lemma10_3}$) and ($\ref{eq:lemma10_4}$),
$$
\begin{aligned}
    \langle \nabla \varphi(x^{++}) - \nabla \varphi(x^{+}), x^{++}-x^{+} - d  \rangle &\leq  0 \\
    \langle \nabla \varphi(x^{++}) - \nabla \varphi(x^{+}), x^{++}-x^{+} \rangle &\leq     \langle \nabla \varphi(x^{++}) - \nabla \varphi(x^{+}), y-x \rangle 
\end{aligned}
$$
By the strong convexity of $\varphi$,
$$
\begin{aligned}
    \mu \|x^{++} - x^{+}\|^{2} &\leq   \langle \nabla \varphi(x^{++}) - \nabla \varphi(x^{+}), y-x \rangle \\ 
    \|x^{++} - x^{+}\|^{2} &\leq  \mu^{-1} \| \nabla \varphi(x^{++}) - \nabla \varphi(x^{+}) \|_{*} \cdot \|x-y\|
\end{aligned}
$$
By the smoothness of $\varphi$,
\begin{subequations}
\label{eq:lemma10_euclidean_proj_expansiveness}
\begin{align}
    \|x^{++} - x^{+}\|^{2} &\leq  (\mu/L)^{-1} \| x^{++} - x^{+} \| \cdot \|x-y\|  \\
    \|x^{++} - x^{+}\|^{2} &\leq  (\mu/L)^{-2}  \|x-y\|^{2} 
\end{align}
\end{subequations}
By the definition of smoothness and strong convexity,
$$
D_{\varphi}(x^{++}, x^{+}) \leq \frac{L}{2}\|x^{++} - x^{+}\|^{2} \text{\quad and \quad }  \frac{\mu}{2} \|x - y\|^{2}\leq D_{\varphi}(x,y),
$$
Together with (\ref{eq:lemma10_euclidean_proj_expansiveness}), 
\begin{equation}
    D_{\varphi}(x^{++}, x^{+}) \leq (L/\mu)^{3}  D_{\varphi}(x, y).
\end{equation}

\hfill \Halmos

\end{proof}

\section{Proof of Proposition \ref{prop: monotone_BPP_is_PPM}}
\label{sec: proof prop: monotone_BPP_is_PPM}
We begin by proving (C1) is a sufficient condition for two paths to coincide: $\{x_{k}\}$ is monotone on $\mathcal{X}$ $\Rightarrow$ $x_{k+1} = x_{\mathrm{{\CentralPathAbb}}}\left(\omega_{k+1}\right), \ \forall k\in [0,K], $
where  $\omega_{k} = \left( \sum_{j=0}^{k-1} \lambda_{j}^{-1} \right)^{-1}$. We proceed with a proof by induction: For $k=0$, by Definition \ref{def:proximal-path} and \ref{def:central-path}, we have $x_{1} = x_{\mathrm{{\CentralPathAbb}}}\left(\lambda_{0}\right)$. Next, we prove $x_{k} = x_{\mathrm{{\CentralPathAbb}}}(\omega_{k})$ and (C1) $\Rightarrow$ $x_{k+1} = x_{\mathrm{{\CentralPathAbb}}}(\omega_{k+1})$.  By Theorem 18.2. of \cite{Rockafellar_Convex_Analysis}, for any $x\in \mathcal{X}$, $\exists \mathcal{F}$ which is a face of $\mathcal{X}$ such that $x \in \mathrm{ri}(\mathcal{F})$. Under this guarantee, let $x_{k} \in \mathrm{ri}(\mathcal{F}')$, where $\mathcal{F}'$ is a face of $\mathcal{X}$. Together with $x_{k} = x_{\mathrm{{\CentralPathAbb}}}\left(\omega_{k}\right)$:
$\left \langle a_{0} + \omega_{k} \left( \nabla \varphi(x_{k}) - \nabla \varphi (x_{0}) \right) , x - x_{k} \right \rangle \geq 0, \ \forall x \in \mathcal{X}, $
we have: $\left \langle a_{0} + \lambda_{0} \left( \nabla \varphi(x_{k}) - \nabla \varphi (x_{0}) \right) , x - x_{k} \right \rangle =  0, \ \forall x \in \mathcal{F}'.$
By the monotonicity of $\{x_{k}\}$ on $\mathcal{X}$, $x_{k+1}\in \mathcal{F}'$, therefore: $\left \langle a_{0} + \lambda_{0} \left( \nabla \varphi(x_{k}) - \nabla \varphi (x_{0}) \right) , x_{k+1} - x_{k} \right \rangle =  0,$
by Claim \ref{claim: ppm_is_bpp}, we have $x_{k+1} = x_{\mathrm{{\CentralPathAbb}}}(\omega_{k+1})$, where  $\omega_{k} = \left( \sum_{j=0}^{k-1} \lambda_{j}^{-1} \right)^{-1}.$

For C2 and C3, it is trivial that if $\mathcal{X}$ is an affine subspace $\{x \in \mathbb{R}^{n}: Ax = b \}$ or the entire vector space $\mathbb{R}^{n}$, then $\{x_{k}\}$ must be monotone on $\mathcal{X}$. 

\begin{claim}
\label{claim: ppm_is_bpp}
 Assume $\mathcal{V}$ satisfies Assumption \ref{assumption: blanket V}, and $\varphi$ satisfies Assumption \ref{assumption: blanket varphi}. Let $\{x_k\}$ be a {\ProximalPathFull} initialized by $x_0$ associated with the step-size sequence $\{\lambda_{k}\}$. Let $\{x_{\mathrm{{\CentralPathAbb}}}(\omega_{k})\}$ be the {\CentralPathFull} initialized also at $x_0$. If $\{x_k\}$ satisfies
 $\left\langle a_{0} + \omega_{k} \left(\nabla \varphi(x_k) - \nabla \varphi(x_0) \right) , x_{k+1} - x_k \right\rangle = 0, \ \forall k\in [1,K], $
 then $x_{k+1} = x_{\mathrm{{\CentralPathAbb}}}\left(\omega_{k+1}\right), \ \forall k\in [0,K],$
 where  $\omega_{k} = \left( \sum_{j=0}^{k-1} \lambda_{j}^{-1} \right)^{-1}.$
\end{claim}
\begin{proof}{Proof.}
For $k=0$, we show $x_{1} = x_{\mathrm{{\CentralPathAbb}}}\left(\lambda_{0}\right)$.
By the variational inequality (V.I.) definition of $x_{1}$ as the proximal point update from $x_0$, we have $$\left \langle a_{0} + \lambda_{0} \left( \nabla \varphi(x_{1}) - \nabla \varphi (x_{0}) \right) , x - x_{1} \right \rangle \geq 0, \ \forall x \in \mathcal{X}, $$
which is also precisely the V.I. definition of for $x_{\mathrm{{\CentralPathAbb}}}\left(\lambda_{0}\right)$, therefore $x_{1} = x_{\mathrm{{\CentralPathAbb}}}\left(\lambda_{0}\right)$.

For $k\in [1, K]$, we provide a proof by induction. We begin by proving for $k=1$, if $\left\langle a_{0} + \lambda_{0} \left(\nabla \varphi(x_1) - \nabla \varphi(x_0) \right) , x_{2} - x_1 \right\rangle = 0$, then $x_2 = x_{\mathrm{{\CentralPathAbb}}}((\lambda_{0}^{-1} + \lambda_{1}^{-1})^{-1})$: Similar to the V.I. definition of $x_{1}$ as the proximal point update from $x_0$,  $\left \langle a_{0} + \lambda_{0} \left( \nabla \varphi(x_{1}) - \nabla \varphi (x_{0}) \right) , x - x_{1} \right \rangle \geq 0, \ \forall x \in \mathcal{X}, $
given $x_{2}$ is the proximal point update from $x_{1}$, we have
$$\left \langle a_{0} + \lambda_{1} \left( \nabla \varphi(x_{2}) - \nabla \varphi (x_{1}) \right) , x - x_{2} \right \rangle \geq 0, \ \forall x \in \mathcal{X}.$$
Or equivalently,
$$ \lambda_{0}^{-1}\langle a_{0}, x \rangle + \langle  \nabla \varphi(x_{1}), x \rangle - \langle \nabla \varphi(x_{0}) , x\rangle - \lambda_{0}^{-1} \langle a_{0} , x_{1} \rangle - \langle \nabla \varphi(x_{1}), x_{1}\rangle + \langle \nabla \varphi(x_{0}), x_{1}\rangle\geq 0, \ \forall x \in \mathcal{X},$$ and 
$$ \lambda_{1}^{-1}\langle a_{0}, x \rangle + \langle  \nabla \varphi(x_{2}), x \rangle - \langle \nabla \varphi(x_{1}) , x\rangle - \lambda_{1}^{-1} \langle a_{0} , x_{2} \rangle - \langle \nabla \varphi(x_{2}), x_{2}\rangle + \langle \nabla \varphi(x_{1}), x_{2}\rangle\geq 0, \ \forall x \in \mathcal{X}.$$
Combining the two inequalities, we have 
$$
\begin{aligned}
&(\lambda_{0}^{-1} +  \lambda_{1}^{-1})\langle a_{0}, x \rangle +  \langle  \nabla \varphi(x_{2}), x \rangle - \langle \nabla \varphi(x_{0}) , x\rangle - \langle \nabla \varphi(x_{2}), x_{2}\rangle  - \lambda_{1}^{-1} \langle a_{0} , x_{2} \rangle \\
&- \lambda_{0}^{-1} \langle a_{0} , x_{1} \rangle - \langle \nabla \varphi(x_{1}), x_{1}\rangle + \langle \nabla \varphi(x_{0}), x_{1}\rangle + \langle \nabla \varphi(x_{1}), x_{2}\rangle \geq 0, \ \forall x \in \mathcal{X}.
\end{aligned}
$$ 
Then, if the following equality holds: $\left\langle a_{0} + \lambda_{0} \left(\nabla \varphi(x_1) - \nabla \varphi(x_0) \right) , x_{2} - x_1 \right\rangle = 0,$
or equivalently
$$- \lambda_{0}^{-1} \langle a_{0} , x_{1} \rangle - \langle \nabla \varphi(x_{1}), x_{1}\rangle + \langle \nabla \varphi(x_{0}), x_{1}\rangle + \langle \nabla \varphi(x_{1}), x_{2}\rangle = - \lambda_{0}^{-1} \langle a_{0} , x_{2} \rangle + \langle \nabla \varphi(x_{0}), x_{2} \rangle,$$
we have 
$$
(\lambda_{0}^{-1} +  \lambda_{1}^{-1})\langle a_{0}, x \rangle +  \langle  \nabla \varphi(x_{2}), x \rangle - \langle \nabla \varphi(x_{0}) , x\rangle - \langle \nabla \varphi(x_{2}), x_{2}\rangle - (\lambda_{0}^{-1} + \lambda_{1}^{-1}) \langle a_{0} , x_{2} \rangle + \langle \nabla \varphi(x_{0}), x_{2} \rangle \geq 0, \ \forall x \in \mathcal{X},
$$
which simplifies to
$$\left\langle a_{0} + (\lambda_{0}^{-1} + \lambda_{1}^{-1})^{-1}(\nabla \varphi(x_{2}) - \nabla \varphi(x_{0}) ), x - x_{2} \right \rangle  \geq 0, \ \forall x \in \mathcal{X},$$
By definition, $x_2 = x_{\mathrm{{\CentralPathAbb}}}((\lambda_{0}^{-1} + \lambda_{1}^{-1})^{-1})$.

To finish the induction, we prove for any $k \in [2,K]$, if the following equality holds, $\left\langle a_{0} + \omega_{k} \left(\nabla \varphi(x_k) - \nabla \varphi(x_0) \right) , x_{k+1} - x_k \right\rangle = 0$ and assume $x_{k} = x_{\mathrm{{\CentralPathAbb}}}\left(\omega_{k}\right)$, then we have $x_{k+1} = x_{\mathrm{{\CentralPathAbb}}}\left(\omega_{k+1}\right)$:
By the V.I. definition of $x_{k} = x_{\mathrm{{\CentralPathAbb}}}\left(\omega_{k}\right)$: 
$$\left \langle a_{0} + \omega_{k} \left( \nabla \varphi(x_{k}) - \nabla \varphi (x_{0}) \right) , x - x_{k} \right \rangle \geq 0, \ \forall x \in \mathcal{X}, $$
given $x_{k+1}$ is the proximal point update from $x_{k}$, we have
$$\left \langle a_{0} + \lambda_{k} \left( \nabla \varphi(x_{k+1}) - \nabla \varphi (x_{k}) \right) , x - x_{k+1} \right \rangle \geq 0, \ \forall x \in \mathcal{X}.$$
Or equivalently,
$$ \omega_{k}^{-1}\langle a_{0}, x \rangle + \langle  \nabla \varphi(x_{k}), x \rangle - \langle \nabla \varphi(x_{0}) , x\rangle - \omega_{k}^{-1} \langle a_{0} , x_{k} \rangle - \langle \nabla \varphi(x_{k}), x_{k}\rangle + \langle \nabla \varphi(x_{0}), x_{k}\rangle\geq 0, \ \forall x \in \mathcal{X},$$ and
$$ \lambda_{k}^{-1}\langle a_{0}, x \rangle + \langle  \nabla \varphi(x_{k+1}), x \rangle - \langle \nabla \varphi(x_{k}) , x\rangle - \lambda_{k}^{-1} \langle a_{0} , x_{k+1} \rangle - \langle \nabla \varphi(x_{k+1}), x_{k+1}\rangle + \langle \nabla \varphi(x_{k}), x_{k+1}\rangle\geq 0, \ \forall x \in \mathcal{X}.$$
Combining the two inequalities, we have 
$$
\begin{aligned}
&(\omega_{k}^{-1} +  \lambda_{k}^{-1})\langle a_{0}, x \rangle +  \langle  \nabla \varphi(x_{k+1}), x \rangle - \langle \nabla \varphi(x_{0}) , x\rangle - \langle \nabla \varphi(x_{k+1}), x_{k+1}\rangle  -\lambda_{k}^{-1} \langle a_{0} , x_{k+1} \rangle \\ 
&- \omega_{k}^{-1} \langle a_{0} , x_{k} \rangle - \langle \nabla \varphi(x_{k}), x_{k}\rangle + \langle \nabla \varphi(x_{0}), x_{k}\rangle + \langle \nabla \varphi(x_{k}), x_{k+1}\rangle \geq 0, \ \forall x \in \mathcal{X},
\end{aligned}
$$
then, if the following equality holds: $$\left\langle a_{0} + \omega_{k} \left(\nabla \varphi(x_k) - \nabla \varphi(x_0) \right) , x_{k+1} - x_k \right\rangle = 0,$$
or equivalently
$$- \omega_{k}^{-1} \langle a_{0} , x_{k} \rangle - \langle \nabla \varphi(x_{k}), x_{k}\rangle + \langle \nabla \varphi(x_{0}), x_{k}\rangle + \langle \nabla \varphi(x_{k}), x_{k+1}\rangle = - \omega_{k}^{-1} \langle a_{0} , x_{k+1} \rangle + \langle \nabla \varphi(x_{0}), x_{k+1} \rangle,$$
we have 
$$
\begin{aligned}
 &(\omega_{k}^{-1} +  \lambda_{k}^{-1})\langle a_{0}, x \rangle +  \langle  \nabla \varphi(x_{k+1}), x \rangle - \langle \nabla \varphi(x_{0}) , x\rangle - \langle \nabla \varphi(x_{k+1}), x_{k+1}\rangle  \\
&- (\omega_{k}^{-1} +  \lambda_{k}^{-1}) \langle a_{0} , x_{k+1} \rangle + \langle \nabla \varphi(x_{0}), x_{k+1} \rangle \geq 0, \ \forall x \in \mathcal{X},   
\end{aligned}
$$
which simplifies to
$$\left\langle a_{0} + (\omega_{k}^{-1} +  \lambda_{k}^{-1})^{-1}(\nabla \varphi(x_{k+1}) - \nabla \varphi(x_{0}) ), x - x_{k+1} \right \rangle  \geq 0, \ \forall x \in \mathcal{X},$$
By definition, $x_{k+1} = x_{\mathrm{{\CentralPathAbb}}}((\omega_{k}^{-1} +  \lambda_{k}^{-1})^{-1}) =x_{\mathrm{{\CentralPathAbb}}}(\omega_{k+1}) $. \hfill \Halmos
\end{proof}

\section{Proof of Proposition \ref{prop: X_matches_V}}
\label{sec: proof prop: X_matches_V}
    Under Assumption \ref{assumption: blanket varphi}, $g$ is essentially strict convex with $\nabla g(0)=0$ $\Rightarrow$ $g$ is monotonically increasing over $\mathbb{R}_{+}$, as a result $\mathcal{X}=\{ x\in \mathbb{R}^{n}: \ \|x\|_{\mathcal{U}} \leq l\}=\{ x\in \mathbb{R}^{n}: \ \varphi(x) \leq g(l)\}$, where $\varphi(x)= g\circ \|x\|_{\mathcal{U}}$. In addition, given the uncertainty set is designed to be $\mathcal{U}^{\circ}$, $\varphi(x)= g\circ \|x\|_{\mathcal{U}}= g\circ \|x\|_{\mathcal{U}^{\circ\circ}}$ is also the distance-generating function inducing the optimization path approximations of the robust path, namely $\{x_k\}$ and $\{x_{\mathrm{{\CentralPathAbb}}}(\omega)\}$. 
    By Definition \ref{def:central-path}, for any $\omega \geq 0$, 
    $x_{\mathrm{{\CentralPathAbb}}}(\omega) = \argmin_{ x\in \mathbb{R}^{n}} \{ \langle a_0,x  \rangle + \omega D_{\varphi}(x,x_{0}): \ \varphi(x) \leq g(l) \}.$
    Applying the KKT condition, we obtain
    $$
    x_{\mathrm{{\CentralPathAbb}}}(\omega) = 
    \left\{ \begin{array}{lcl}
    \nabla \varphi^{*}\left( \nabla \varphi(x_{0}) - \omega^{-1} a_{0} \right) & \quad  \text{if} \ \varphi(x_{\mathrm{{\CentralPathAbb}}}(\omega)) < g(l) \\ 
    \nabla  \varphi^{*}\left( \theta_{\omega} \left(\nabla \varphi(x_{0}) - \omega^{-1} a_{0} \right) \right) & \quad \text{if} \ \varphi(x_{\mathrm{{\CentralPathAbb}}}(\omega)) = g(l) \\
    \end{array} \right. ,
    $$
    where $\theta_{\omega} = \left(1+\mu \omega^{-1}\right)^{-1} $ for some Lagrangian dual variable $\mu = \mu(\omega)>0$ such that $\varphi(x_{\mathrm{{\CentralPathAbb}}}(\omega)) = g(l) $. Geometrically, $\theta_{\omega}$ can be interpreted as a rescaling factor of the dual space vector $\left(\nabla \varphi(x_{0}) - \omega^{-1} a_{0}\right)$ such that $\varphi(x_{\mathrm{{\CentralPathAbb}}}(\omega)) = g(l)$, leading to $x_{\mathrm{{\CentralPathAbb}}}(\omega) = \Pi^{\varphi}_{\{ x\in \mathbb{R}^{n}: \ \varphi(x) \leq g(l)\}} \left(\nabla\varphi^{*}\left( \theta_{\omega} \left(\nabla \varphi(x_{0}) - \omega^{-1} a_{0} \right) \right) \right) =\nabla\varphi^{*}\left( \theta_{\omega} \left(\nabla \varphi(x_{0}) - \omega^{-1} a_{0} \right) \right) $.

    Apply a similar procedure for the $\{x_{k}\}$ yields for any $k$, 
    $$
    x_{k} = 
    \left\{ \begin{array}{lcl}
    \nabla \varphi^{*}\left( \nabla \varphi(x_{k-1}) - \lambda_{k}^{-1} a_{0} \right) & \quad  \text{if} \ \varphi(x_{k}) < g(l) \\ 
    \nabla  \varphi^{*}\left( \theta_{\lambda_{k}} \left(\nabla \varphi(x_{k-1}) - \lambda_{k}^{-1} a_{0} \right) \right) & \quad \text{if} \ \varphi(x_{k}) = g(l) \\
    \end{array} \right. ,
    $$
    where $\theta_{\lambda_{k}} = \left(1+\mu \lambda_{k}^{-1}\right)^{-1}$ for some Lagrangian dual variable $\mu = \mu(\lambda_{k})>0$ such that $\varphi(x_k) = g(l) $. 
    
    Given the above characterization of $\{x_{\mathrm{{\CentralPathAbb}}}(\omega)\}$ and $\{x_{k}\}$, we now prove $\{x_k\} \subset \{x_{\mathrm{{\CentralPathAbb}}}\left(\omega\right): \omega \in [0,\infty)\}$ by showing for any $x_k \in \{x_k\}$, there exists a $x_{\mathrm{{\CentralPathAbb}}}\left(\omega_{k}\right) \in \{x_{\mathrm{{\CentralPathAbb}}}\left(\omega\right): \omega \in [0,\infty)\}$ such that $x_k = x_{\mathrm{{\CentralPathAbb}}}\left(\omega_{k}\right)$. We proceed with a proof by induction. 
    
    For $k=1$: by Definition \ref{def:proximal-path} and \ref{def:central-path}, we have $x_1 = x_{\mathrm{{\CentralPathAbb}}}(\lambda_0)$. 
    
    For $k>1$: assume $x_k =x_{\mathrm{{\CentralPathAbb}}}(\omega_k)$, in the dual space, w.l.o.g., assume $\varphi(x_k) < g(l)$, we have 
    $$
    \nabla \varphi(x_k) = \nabla \varphi(x_{0}) - \omega_{k}^{-1} a_{0}.  
    $$ Next, given $x_{k+1}$ as the proximal point update of $x_k$, w.l.o.g., assume $\varphi(x_{k+1}) = g(l)$ yields
    $$
    \begin{aligned}
    \nabla \varphi(x_{k+1}) &= 
     \theta_{\lambda_{k+1}} \left(\nabla \varphi(x_{k}) - \lambda_{k+1}^{-1} a_{0} \right)  \\
     &= 
      \theta_{\lambda_{k+1}} \left(    \left( \nabla \varphi(x_{0}) - \omega_{k}^{-1} a_{0} \right) - \lambda_{k+1}^{-1} a_{0} \right)  \\
      &= 
      \theta_{\lambda_{k+1}} \left( \nabla \varphi(x_{0}) - \left(\omega_{k}^{-1} + \lambda_{k+1}^{-1} \right) a_{0}  \right) \\ 
    &=\nabla \varphi\left( x_{\mathrm{{\CentralPathAbb}}}\left(\omega_{k}^{-1} + \lambda_{k+1}^{-1} \right) \right).   
    \end{aligned}
    $$
By invoking again Lemma \ref{lemma: legendre_bijective} and mapping both the LHS and RHS back to the primal space, we have $x_{k+1} \in \{x_{\mathrm{{\CentralPathAbb}}}\left(\omega\right): \omega \in [0,\infty)\}$. 

\section{Proof of Theorem \ref{theorem: PPM as approx BPP}}
\label{sec: proof theorem: PPM as approx BPP}
We begin by constructing for each $i\in [I]$, two auxiliary {\CentralPathFull}s that mirror exactly $\left\{x_k: k\in \left[\underline{k}^{(i)}+1, \overline{k}^{(i)}\right]\right\}$ and $\left\{x_{\mathrm{{\CentralPathAbb}}}\left(\upsilon_{k}^{-1}; x_{0} \right): k\in \left[\underline{k}^{(i)}+1, \overline{k}^{(i)}\right]\right\}$ respectively. The distance between the two original sequences can be analyzed equivalently via the distance between the two auxiliary sequences.

\textbf{Step One:} We first construct the auxiliary {\CentralPathFull}  for $\left\{x_k: k\in \left[\underline{k}^{(i)}+1, \overline{k}^{(i)}\right]\right\}$. By Definition \ref{def:monotone}, $\left\{x_k: k\in \left[\underline{k}^{(i)}+1, \overline{k}^{(i)}\right]\right\}$ is monotone on $\mathrm{ri}(\mathcal{F}_{i})$, consequently by Proposition \ref{prop: monotone_BPP_is_PPM} we have
$$x_k = x_{\mathrm{{\CentralPathAbb}}}\left(\left(\upsilon_{k}' \right)^{-1}; x_{\underline{k}^{(i)}}\right), \quad \forall k\in \left[\underline{k}^{(i)}+1, \overline{k}^{(i)}\right],$$
where $\upsilon_{k}'  =  \sum_{j =\underline{k}^{(i)}}^{k-1} \lambda_{j}
^{-1} $. It is helpful to consider its Bregman projection reformulation due to Theorem \ref{theorem: RP_CP_PP are Bregman Projections}: $x_{\mathrm{{\CentralPathAbb}}}\left(\left(\upsilon_{k}' \right)^{-1}; x_{\underline{k}^{(i)}}\right) = \Pi^{\varphi}_{\mathrm{ri}(\mathcal{F}_{i})}\left( \nabla \varphi^{*}\left(\nabla \varphi\left(x_{\underline{k}^{(i)}}\right) -\upsilon_{k}'a_{0}\right)\right)$.

\ 

\textbf{Step Two:} Next we construct the auxiliary {{\CentralPathFull}} for $\left\{x_{\mathrm{{\CentralPathAbb}}}\left(\upsilon_{k}^{-1}; x_{0} \right): k\in \left[\underline{k}^{(i)}+1, \overline{k}^{(i)}\right]\right\}$ where $\upsilon_{k} = \underline{\upsilon}^{(i)} +  \sum_{j =\underline{k}^{(i)}}^{k-1} \lambda_{j}^{-1}$, and show the two sequences are equivalent. 
First, we define the following point that is critical:
$$x_{\mathrm{{\CentralPathAbb}}}\left((\underline{\upsilon}^{(i)})^{-1}; x_0 \right) = \Pi^{\varphi}_{\mathrm{ri}(\mathcal{F}_{i})}( \underbrace{\nabla \varphi^{*}(\nabla \varphi\left(x_{0}\right) -\underline{\upsilon}^{(i)}a_{0})}_{\bigstar}).$$

By Theorem \ref{theorem: RP_CP_PP are Bregman Projections} and given $\upsilon_{k} = \underline{\upsilon}^{(i)} +  \sum_{j =\underline{k}^{(i)}}^{k-1} \lambda_{j}^{-1}$, $\left\{x_{\mathrm{{\CentralPathAbb}}}\left(\upsilon_{k}^{-1}; x_{0} \right): k\in \left[\underline{k}^{(i)}+1, \overline{k}^{(i)}\right]\right\}$ has the following Bregman projection reformulation,
$$x_{\mathrm{{\CentralPathAbb}}}\left(\upsilon_{k}^{-1}; x_{0} \right) = \Pi^{\varphi}_{\mathrm{ri}(\mathcal{F}_{i})}\left( \nabla \varphi^{*} \left( \underbrace{(\nabla \varphi(x_{0}) -\underline{\upsilon}^{(i)}a_{0})}_{\nabla \varphi(\bigstar)} - \left(\sum_{j =\underline{k}^{(i)}}^{k-1} \lambda_{j}^{-1}\right) a_{0}\right) \right).$$
We define its auxiliary sequence as
$$x_{\mathrm{{\CentralPathAbb}}}\left((\upsilon_{k}' )^{-1}; x_{\mathrm{{\CentralPathAbb}}}\left((\underline{\upsilon}^{(i)})^{-1}; x_0 \right) \right) = \Pi^{\varphi}_{\mathrm{ri}(\mathcal{F}_{i})}\left( \nabla \varphi^{*} \left( \underbrace{\nabla \varphi\left(x_{\mathrm{{\CentralPathAbb}}}\left((\underline{\upsilon}^{(i)})^{-1}; x_0 \right)\right)}_{\nabla \varphi\left(\Pi^{\varphi}_{\mathrm{ri}(\mathcal{F}_{i})}\left(\bigstar\right)\right)} - \left(\sum_{j =\underline{k}^{(i)}}^{k-1} \lambda_{j}^{-1}\right) a_{0}\right) \right),$$
where $\upsilon_{k}'  =  \sum_{j =\underline{k}^{(i)}}^{k-1} \lambda_{j}
^{-1} $.

By definition, $\Pi^{\varphi}_{\mathcal{X}}\left( \bigstar \right) \in \mathrm{ri}(\mathcal{F}_{i})$, therefore by Claim \ref{claim: proj_on_X_in_riF} we have
$\Pi^{\varphi}_{\mathrm{ri}(\mathcal{F}_i)}(\bigstar) = \Pi^{\varphi}_{\mathrm{Aff}(\mathrm{ri}(\mathcal{F}_i))}(\bigstar).$ Consequently, by the same argument as in the proof of Corollary \ref{corollary: BPP is RP}, two sequences are equivalent, and we have
$$x_{\mathrm{{\CentralPathAbb}}}\left(\upsilon_{k}^{-1}; x_{0} \right) = x_{\mathrm{{\CentralPathAbb}}}\left((\upsilon_{k}' )^{-1}; x_{\mathrm{{\CentralPathAbb}}}\left((\underline{\upsilon}^{(i)})^{-1}; x_0 \right) \right), \quad \forall k\in \left[\underline{k}^{(i)}+1, \overline{k}^{(i)}\right],$$
where $\upsilon_{k} = \underline{\upsilon}^{(i)} +  \sum_{j =\underline{k}^{(i)}}^{k-1} \lambda_{j}^{-1}$ and $\upsilon_{k}' =   \sum_{j =\underline{k}^{(i)}}^{k-1} \lambda_{j}^{-1}$.

\textbf{Step Three:} By the results established in the previous two steps, for each $k\in \left[\underline{k}^{(i)}+1, \overline{k}^{(i)}\right]$ we have
$$
\begin{aligned}
    &D_{\varphi} \left(x_{k}, x_{\mathrm{{\CentralPathAbb}}}\left(\upsilon_k^{-1}; x_{0}\right) \right) \\
    = &D_{\varphi} \left(x_{\mathrm{{\CentralPathAbb}}}\left((\upsilon_{k}' )^{-1}; x_{\underline{k}^{(i)}}\right), x_{\mathrm{{\CentralPathAbb}}}\left((\upsilon_{k}' )^{-1}; x_{\mathrm{{\CentralPathAbb}}}\left((\underline{\upsilon}^{(i)})^{-1}; x_0 \right) \right) \right)\\
   = &D_{\varphi} \left(\Pi^{\varphi}_{\mathrm{ri}(\mathcal{F}_{i})}\left( \nabla \varphi^{*}\left(\nabla \varphi\left(x_{\underline{k}^{(i)}}\right) -\upsilon_{k}'a_{0}\right)\right) , \Pi^{\varphi}_{\mathrm{ri}(\mathcal{F}_{i})}\left( \nabla \varphi^{*} \left( \nabla \varphi\left(x_{\mathrm{{\CentralPathAbb}}}\left((\underline{\upsilon}^{(i)})^{-1}; x_0 \right)\right) - \upsilon_{k}' a_{0}\right) \right)\right) \\
   = &D_{\varphi} \left(\Pi^{\varphi}_{\mathrm{Aff}(\mathrm{ri}(\mathcal{F}_{i}))}\left( \nabla \varphi^{*}\left(\nabla \varphi\left(x_{\underline{k}^{(i)}}\right) -\upsilon_{k}'a_{0}\right)\right) , \Pi^{\varphi}_{\mathrm{Aff}(\mathrm{ri}(\mathcal{F}_{i}))}\left( \nabla \varphi^{*} \left( \nabla \varphi\left(x_{\mathrm{{\CentralPathAbb}}}\left((\underline{\upsilon}^{(i)})^{-1}; x_0 \right)\right) - \upsilon_{k}' a_{0}\right) \right)\right) \\
   \leq & \kappa \cdot  D_{\varphi}\left(x_{\mathrm{{\CentralPathAbb}}}\left((\underline{\upsilon}^{(i)})^{-1}; x_0 \right), x_{\underline{k}^{(i)}} \right). 
\end{aligned} 
$$ 
where the third equality is due to Claim \ref{claim: proj_on_X_in_riF} and the inequality is due to Claim \ref{claim: Breg_Proj_calc_2}. \hfill \Halmos

\begin{claim}
\label{claim: proj_on_X_in_riF}
Let $\mathcal{X}$ be a closed convex polyhedron, and let $\mathcal{F}$ be a face of $\mathcal{X}$. Additionally, assume $\varphi$ is Legendre, then
$$ y = \Pi^{\varphi}_{\mathcal{X}}(x) \in \mathrm{ri}(\mathcal{F}) \quad \Rightarrow \quad y = \Pi^{\varphi}_{\mathrm{ri}(\mathcal{F})}(x) = \Pi^{\varphi}_{\mathrm{Aff}(\mathrm{ri}(\mathcal{F}))}(x).$$
\end{claim}
\begin{proof}{Proof of Claim \ref{claim: proj_on_X_in_riF}.}
Let $\mathcal{N}_{\mathcal{C}}(p)$ denote the normal cone of a convex set $\mathcal{C}$ at $p\in \mathcal{C}$, by Lemma \ref{lemma: Bregman Projection VI}, we have
$$\nabla \varphi(y) - \nabla \varphi(x) \in \mathcal{N}_{\mathcal{X}}(y).$$
The result follows by proving the following statement is true:
$$\mathcal{N}_{\mathcal{X}}(y) \subset \mathcal{N}_{\mathrm{ri}(\mathcal{F})}(y) = \mathcal{N}_{\mathrm{Aff}(\mathrm{ri}(\mathcal{F}))}(y).$$
We begin by showing $\mathcal{N}_{\mathcal{X}}(y) \subset \mathcal{N}_{\mathrm{ri}(\mathcal{F})}(y)$. Assume $d \in \mathcal{N}_{\mathcal{X}}(y)$, by definition we have 
$$\langle d, y-x \rangle, \quad \forall x\in \mathcal{X},$$
given $\mathrm{ri}(\mathcal{F}) \subset \mathcal{X}$ we have
$$\langle d, y-x \rangle, \quad \forall x\in \mathrm{ri}(\mathcal{F}),$$
hence $d \in \mathcal{N}_{\mathrm{ri}(\mathcal{F})}(y)$.

By the same argument it follows that $\mathcal{N}_{\mathrm{Aff}(\mathrm{ri}(\mathcal{F}))}(y) \subset \mathcal{N}_{\mathrm{ri}(\mathcal{F})}(y)$. It remains to show $\mathcal{N}_{\mathrm{ri}(\mathcal{F})}(y) \subset \mathcal{N}_{\mathrm{Aff}(\mathrm{ri}(\mathcal{F}))}(y)$:
Assume $g\in \mathcal{N}_{\mathrm{ri}(\mathcal{F})}(y)$, by definition we have
$$\langle g ,z-y \rangle \leq 0, \quad \forall z \in \mathrm{ri}(\mathcal{F}).$$
For any $w \in \mathrm{Aff}(\mathrm{ri}(\mathcal{F}))$, , there exist some $\alpha \geq 0$ such that
$$y + \alpha(w-y) \in \mathrm{ri}(\mathcal{F}),$$
consequently,
$$
\begin{aligned}
    \langle g, y + \alpha(w-y) - y &\rangle \leq 0, \quad \quad (g\in \mathcal{N}_{\mathrm{ri}(\mathcal{F})}(y))\\
    \alpha \langle g,w-y &\rangle \leq 0, \\
    \langle g, w-y &\rangle \leq 0.  \quad \quad  (\alpha \geq 0)\\
\end{aligned}
$$
hence $g \in \mathcal{N}_{\mathrm{Aff}(\mathrm{ri}(\mathcal{F}))}(y)$.
We have shown,
$$\mathcal{N}_{\mathcal{X}}(y) \subset \mathcal{N}_{\mathrm{ri}(\mathcal{F})}(y) = \mathcal{N}_{\mathrm{Aff}(\mathrm{ri}(\mathcal{F}))}(y).$$
To finish the proof, given
$$\nabla \varphi(y) - \nabla \varphi(x) \in \mathcal{N}_{\mathcal{X}}(y),$$
by $\mathcal{N}_{\mathcal{X}}(y) \subset \mathcal{N}_{\mathrm{ri}(\mathcal{F})}(y) = \mathcal{N}_{\mathrm{Aff}(\mathrm{ri}(\mathcal{F}))}(y)$, we have 
$$\nabla \varphi(y) - \nabla \varphi(x) \in \mathcal{N}_{\mathrm{ri}(\mathcal{F})}(y)\quad \text{and}\quad \nabla \varphi(y) - \nabla \varphi(x) \in \mathcal{N}_{\mathrm{Aff}(\mathrm{ri}(\mathcal{F}))}.$$
Or equivalently,
$$y = \Pi^{\varphi}_{\mathrm{ri}(\mathcal{F})}(x) = \Pi^{\varphi}_{\mathrm{Aff}(\mathrm{ri}(\mathcal{F}))}(x).$$ \hfill \Halmos
\end{proof}

\begin{claim}
\label{claim: Breg_Proj_calc_2}
Let $\mathcal{A}\subset \mathbb{R}^{n} $ be a affine subspace, let $x,y \in \mathcal{A}$ and let $\varphi$ be a Legendre function whose Bregman projection is $\kappa-\text{expansive}$, then for any $c \in \mathbb{R}^{n}$,
$$D_{\varphi}\left( \Pi^{\varphi}_{\mathcal{A}}\left( \nabla \varphi^{*}\left(\nabla \varphi(x) -\upsilon c\right)\right), \Pi^{\varphi}_{\mathcal{A}}\left( \nabla \varphi^{*}\left(\nabla \varphi(y) -\upsilon c\right)\right) \right) \leq \kappa\cdot  D_{\varphi}(y,x), \quad \forall \upsilon \in [0,\infty). $$
\end{claim}
\begin{proof}{Proof of Claim \ref{claim: Breg_Proj_calc_2}.}
Denote $\mathcal{A} = \mathcal{L} + x$, where $\mathcal{L} = \mathcal{A} -x $ is a linear subspace. For any $\upsilon \in [0,\infty)$, by Lemma \ref{lemma:dual_cone bregman_proj}, we have 
$$
\begin{aligned}
\nabla \varphi\left(\Pi^{\varphi}_{\mathcal{A}}\left( \nabla \varphi^{*}\left(\nabla \varphi(x) -\upsilon c\right)\right)\right) =& \nabla \varphi\left(\Pi^{\varphi}_{\mathcal{L}+x}\left( \nabla \varphi^{*}\left(\nabla \varphi(x) -\upsilon c\right)\right)\right) \\
=& \Pi^{\varphi^{*}}_{\mathcal{L}^{\bot}+\nabla\varphi(x) -\upsilon c}\left( \nabla \varphi(x)\right), \\
\end{aligned}
$$
and 
$$
\begin{aligned}
\nabla \varphi\left(\Pi^{\varphi}_{\mathcal{A}}\left( \nabla \varphi^{*}\left(\nabla \varphi(y) -\upsilon c\right)\right)\right) 
=& \nabla \varphi\left(\Pi^{\varphi}_{\mathcal{L}+x}\left( \nabla \varphi^{*}\left(\nabla \varphi(y) -\upsilon c\right)\right)\right) \\
=& \Pi^{\varphi^{*}}_{\mathcal{L}^{\bot}+\nabla\varphi(y) -\upsilon c}\left( \nabla \varphi(x)\right). \\
\end{aligned}
$$
With the above results, next we map the pair $ \Pi^{\varphi}_{\mathcal{A}}\left( \nabla \varphi^{*}\left(\nabla \varphi(x) -\upsilon c\right)\right)$ and $\Pi^{\varphi}_{\mathcal{A}}\left( \nabla \varphi^{*}\left(\nabla \varphi(y) -\upsilon c\right)\right)$ to the dual space before mapping them back to the primal space. For any $\upsilon \in [0,\infty)$, 
$$
\begin{aligned}
&D_{\varphi}\left( \Pi^{\varphi}_{\mathcal{A}}\left( \nabla \varphi^{*}\left(\nabla \varphi(x) -\upsilon c\right)\right), \Pi^{\varphi}_{\mathcal{A}}\left( \nabla \varphi^{*}\left(\nabla \varphi(y) -\upsilon c\right)\right) \right)  \\
=&D_{\varphi^{*}}\left( \nabla \varphi\left(\Pi^{\varphi}_{\mathcal{A}}\left( \nabla \varphi^{*}\left(\nabla \varphi(y) -\upsilon c\right)\right)\right), \nabla \varphi\left(\Pi^{\varphi}_{\mathcal{A}}\left( \nabla \varphi^{*}\left(\nabla \varphi(x) -\upsilon c\right)\right)\right) \right)  \\
=&D_{\varphi^{*}}\left(\Pi^{\varphi^{*}}_{\mathcal{L}^{\bot}+\nabla\varphi(y) -\upsilon c}\left( \nabla \varphi(x)\right), \Pi^{\varphi^{*}}_{\mathcal{L}^{\bot}+\nabla\varphi(x) -\upsilon c}\left( \nabla \varphi(x)\right)\right) \\
\leq & \kappa \cdot D_{\varphi^{*}}\left(\nabla\varphi(x), \nabla\varphi(y) \right) \\
=&  \kappa \cdot D_{\varphi}\left(y,x \right), \\
\end{aligned}
$$
where the first equality is due to \cite{bauschke1997legendre} Theorem 3.7(v), the inequality is by Definition \ref{def: Kappa_expansive} and the third equality is again due to \cite{bauschke1997legendre} Theorem 3.7(v).
\hfill \Halmos\end{proof}

\end{APPENDICES}


\end{document}